\documentclass[sn-chicago]{sn-jnl}

\usepackage{graphicx}%
\usepackage{multirow}%
\usepackage{amsmath,amssymb,amsfonts}%
\usepackage{amsthm}%
\usepackage{mathrsfs}%
\usepackage[title]{appendix}%
\usepackage{xcolor}%
\usepackage{textcomp}%
\usepackage{manyfoot}%
\usepackage{booktabs}%
\usepackage{algorithm}%
\usepackage{algorithmicx}%
\usepackage{algpseudocode}%
\usepackage{listings}%
\usepackage{paralist}
\usepackage{amsthm}
\usepackage{tikz}
\usepackage{pgfplots}
\usetikzlibrary{calc}
\pgfplotsset{compat=1.18}

\newtheorem{lemma}{Lemma}

\begin{document}

\title[Article Title]{Robust Condition-Based Operations \& Maintenance: 
Synergizing Multi-Asset Degradation Rate Interactions and Operations-Induced Degradation}


\author[1]{\fnm{Deniz} \sur{Altinpulluk}}\email{deniz@wayne.edu}

\author[1]{\fnm{Farnaz} \sur{Fallahi}}\email{farnaz.fallahi@wayne.edu}

\author*[1]{\fnm{Murat} \sur{Yildirim}}\email{murat@wayne.edu}

\author[2]{\fnm{Mohammad Javad} \sur{Feizollahi}}\email{mfeizollahi@gsu.edu}

\affil[1]{\orgdiv{Industrial and Systems Engineering}, \orgname{Wayne State University}, 
\orgaddress{\street{4815 4th St}, \city{Detroit}, \postcode{48201}, \state{MI}, \country{USA}}}

\affil[2]{\orgdiv{Robinson College of Business}, 
\orgname{Georgia State University}, \orgaddress{\street{35 Broad St NW}, \city{Atlanta}, \postcode{30303}, \state{GA}, \country{USA}}}


\abstract{Effective operations and maintenance (O\&M) in modern production systems hinges on careful orchestration of economic and degradation dependencies across a multitude of assets. While the economic dependencies are well studied, degradation dependencies and their impact on system operations remain an open challenge. To address this challenge, we model condition-based production and maintenance decisions for multi-asset systems with degradation interactions. There is a rich literature on condition-based O\&M policies for single-asset systems. These models fail to represent modern systems composed of multiple interacting assets. We are providing the first O\&M model to optimize O\&M in multi-asset systems with embedded decision-dependent degradation interactions. We formulate robust optimization models that inherently capture degradation and failure risks by embedding degradation signals via a set of constraints, and building condition-based uncertainty sets to model probable degradation scenarios. We offer multiple reformulations and a solution algorithm to ensure computational scalability. Performance of the proposed O\&M model is evaluated through extensive experiments, where the degradation is either emulated or taken from vibration-based readings from a rotating machinery system. The proposed model provides significant improvements in terms of operation, maintenance, and reliability metrics. Due to a myriad of dependencies across assets and decisions, it is often difficult to translate \textit{asset-level} failure predictions to \textit{system-level} O\&M decisions. This challenge puts a significant barrier to the return on investment in condition monitoring and smart maintenance systems. Our approach offers a seamless integration of data-driven failure modeling and mathematical programming to bridge the gap across predictive and prescriptive models.}




\maketitle

\section{Introduction}\label{sec1}
Operations and maintenance (O\&M) is a pivotal component of manufacturing and service processes with far-reaching implications for system-level profitability, reliability, and resilience. Owing to its significance, O\&M activities constitute a major market for any industry, with a significant presence in manufacturing, mobility, and energy sectors. The goal of O\&M policies is to effectively use equipment lifetime, prevent unexpected failures, and devise strategies to mitigate downtime or failure-induced operational interruptions. This is a difficult challenge for multi-asset systems due to the existence of many interdependencies across assets and decisions. 
For example, O\&M decisions have significant two-way interactions in most applications. On the one hand, maintenance impact on operations occurs due to interruptions to production and service processes. For instance, the replacement of a \$5,000 bearing in wind farms could cost up to \$250,000 due to operational interruptions  \citep{hu2020engineering}. Likewise, maintenance costs can constitute 15\% to 70\% of total production costs \citep{jafar2022robust}. On the other hand, operations impact maintenance due to additional stress/loading on assets that evolves according to how the components are operated. Asset degradation is significantly impacted by the condition of the other assets and operational decisions in the fleet. For instance, increasing production puts additional stress on certain machines, further accelerating their rate of degradation \citep{uit2020condition}. 

O\&M decisions are still predominantly based on fixed time-based periodic schedules that ignore the operations-maintenance interactions and fail to leverage sensor data—think of the clich\'{e}, change the oil of your car every 3000 miles regardless of car age, make, driving habits, etc. Time-based schedules do not account for the actual condition of the asset/equipment when planning maintenance activities and, therefore, cannot be used to anticipate failures. If implemented in a conservative fashion, time-based policies still drive up the cost of maintenance due to frequent unnecessary maintenances and downtime. To address this issue, there is an increasing interest in the use of condition-based maintenance (CBM) policies, which leverage sensor information to understand asset conditions and O\&M outcomes. A central focus of these policies is to analyze sensor information to discover latent signs of degradation and conduct maintenance as needed. Despite industry interest in condition-based O\&M models, there are significant operations management challenges that hinder their widespread adoption. A significant challenge is that the models used for sensor-driven condition assessment of failure likelihoods and risks are not well-integrated into decision optimization models for multi asset systems. It becomes a central question to determine how the sensor-driven insights on asset conditions and failure risks can be incorporated within scalable and tractable decision models that optimize fleet-level O\&M decisions subject to multiple layers of operational and maintenance dependencies.

Many critical industrial assets are equipped with sensors to monitor their condition and performance ---a prerequisite for the condition-based O\&M models. Raw signals from these sensors often exhibit characteristic features and fault-based patterns that capture information about the physical and performance degradation of these assets. \textit{Condition Monitoring} is the process of collecting sensor data, such as temperature, vibration, noise, etc., from these assets to assess their state of health. For instance, as the rotating components of generators degrade, their vibration levels increase. These vibration readings can be used to assess the level of degradation in these components. These sensor signals are used to generate inferences on asset conditions, which are evaluated through metrics called \textit{degradation signals}. Degradation signals are correlated with degradation severity and provide valuable information on the current state of health and its future trajectory.

Degradation signals are typically used for two categories of predictive models. The first category, called \textit{asset diagnostics}, focuses on estimating the current state of health (or condition) of the asset. The second category, called \textit{asset prognostics}, predicts failure risks and derives the remaining life distribution (RLD) of the assets. Prognostics require both an estimation of the current state and a prediction of future degradation trajectory. Estimation of the degradation trajectory hinges on \textit{degradation modeling} efforts that offer a stochastic model that mimics the degradation trajectory over time. It is often difficult to predict this trajectory, specifically in assets with complex degradation processes that depend on operational decisions and asset-to-asset interactions.

Translating asset prognostics predictions into O\&M models constitutes a significant modeling challenge. Most approaches to CBM rely on policy-based models, such as Markov Decision Models, with limited flexibility to incorporate complex operational constraints and fleet-level interactions, thus focusing mainly on single-asset systems.  \cite{yildirim2017integrated} showcase that these models do not generalize to fleet-level O\&M models, and simpler time-based models that incorporate fleet-level interactions can often outperform the generalizations of CBM models that focus on a single asset at a time. Fleet-level O\&M interactions across assets and decisions can be categorized into two classes: \textit{economic dependencies} and \textit{degradation dependencies}. Economic dependencies relate to a range of interactions that deal with process-driven limitations. For instance, maintenance of a certain asset may require dismantling another asset in the fleet, or operational requirements may dictate that a certain number of assets have to always be available for production. A representative example in offshore wind farm maintenance relates to maintenance logistics, where maintenance of a single turbine requires the transportation of equipment and maintenance crew through work boats. When the crew accesses the location, it creates an economic incentive to maintain multiple turbines ---a concept called opportunistic maintenance. 

The second class of interactions, called degradation dependencies, refers to the interactions that impact how assets degrade and fail; and can occur either due to operational stress or multi-asset interactions. Operational stress refers to cases where the increasing production rate of a machine sparks an increase in the degradation rate of its constituent components. Multi-asset interactions occur when an increase in asset degradation status triggers a further acceleration of degradation in other assets/components in the system. The impact of multi-asset interactions and operational stress on degradation rate manifest in a variety of multi-asset/component systems such as manufacturing equipment, wind turbines, gearboxes, and robotic systems. 
The benefit of considering both impacts in degradation and optimization models is two-fold. First, it provides a more accurate prediction of remaining life distribution that helps reduce the risk of unexpected failures and increases the lifetime utilization of assets. Second and foremost, it allows the decision maker to control the degradation process of an asset by finetuning production rates and fleet-level maintenance schedules. This enables operators to delay or expedite maintenance times as needed.

In this article, we propose an integrated framework for condition-based O\&M at fleet-level, which embeds sensor-driven insights on asset conditions and degradation models within large-scale robust optimization models. This integration provides a direct path to translate sensor-driven insights into fleet-level operational decisions. The proposed approach offers a fundamental shift away from existing O\&M models that either ignore degradation models (e.g., operations-focused time-based maintenance models) or focus on a single asset (e.g., policy-based decision models). Thus, the main contributions of this article can be summarized as follows:
\begin{itemize}
    \item We develop an integrated modeling framework that embeds dynamically updated degradation models within a robust optimization formulation through two modeling innovations. First, we reformulate continuous stochastic functions used for degradation modeling into a set of linear constraints. Second, we incorporate degradation model stochasticity within novel degradation-driven uncertainty sets. The proposed constraints and uncertainty sets adapt to sensor-driven degradation parameters to predict asset conditions and future failure risks.
    \item We explicitly model the significant dependencies across decisions and assets within fleet-level O\&M. In degradation dependencies, we jointly consider operational stress/loading and multi-asset degradation interactions. Incorporating these dependencies require modeling an explicit link across asset-specific degradation constraints and operational decisions. Our approach also enables the consideration of a wide range of operational dependencies.
    \item We formulate a robust optimization problem that captures a degradation-induced nested uncertainty set which allows us to adapt to a wider variety of worst-case realizations and model multiple maintenances in a planning horizon.
    \item We develop extensive computational experiments to demonstrate the value of the proposed decision framework in a range of settings. Our experiments use real-world vibration-based degradation signals from a rotating machinery application to model asset degradation. Degradation signals have been captured from a brand new stage to failure during an accelerated life testing (ALT) procedure conducted under different operational stress/loading environments.
\end{itemize}

\section{Literature Review}
The O\&M framework in this article brings together three streams of research:  
(i) condition monitoring, (ii) degradation modeling, and (iii) condition-based maintenance. 
The primary goal of condition monitoring is to discover and record latent patterns in sensor information that correlates with the asset state of health. Typically, signal processing and pattern recognition techniques are used to process data from mechanical sensors (e.g., vibration analysis to detect cracks), electrical sensors (e.g., partial discharge measurements to evaluate insulation), and chemical sensors (e.g., particle matter concentration in contact points to detect wear and tear). 
The literature on condition monitoring has two main approaches. Data-driven approaches solely rely on statistical and machine learning methods 
that harness historical data to estimate the current health status of components. In contrast, domain-driven approaches build on an in-depth understanding of the underlying component-specific physics-of-failure processes. Domain-driven approaches offer additional predictive stability and enable predictions with limited data yet require a lengthy and expensive development cycle. Condition monitoring has wide applications, including wind turbines \citep{lu2009review}, electric motors \citep{lee2020condition}, and nuclear power plant components \citep{zio2010data}. \cite{tidriri2016bridging} provides a comprehensive review of the subject.

A significant limitation of the condition monitoring systems is their focus on the current state of health. Anticipating future failure risks requires stochastic formulations that can help anticipate how future degradation is likely to evolve. To this end, degradation modeling approaches build stochastic formulations to model the long-term behavior of asset degradation and health \citep{nelson2009accelerated}.  
The stochastic degradation models take the form of either discrete-state (Compound Poisson, Markovian, or semi-Markovian process) or continuous-state (Wiener, Gamma, or Inverse Gaussian process) stochastic processes based on the nature of the degradation  \citep{alaswad2017review}. 
The majority of the degradation modeling approaches assume a constant environment and ignore the significant impact of other components and operational conditions \citep{gebraeel2005residual}. 

In reality, the time-varying conditions due to operational stress or multi-asset interactions can cast a significant impact on the degradation trajectory.

Operational stress refers to dynamic loading on assets due to the mode of operations. For example, manufacturing stations work under different production rates and exhibit various degradation behaviors  
\citep{uit2020condition}. 
Different methodologies have been used to model the effect of time-varying operational stress by incorporating additional parameters and capturing the changes through randomly occurring shocks, i.e., incorporating the impact of prevailing stochastic conditions by modeling them as shocks 
\citep{nakagawa2007shock} and wear process models that consider continuous damage accumulation over time \citep{bian2015degradation,gebraeel2008prognostic}.

Multi-asset degradation interactions constitute the second major factor that impacts the pace of degradation. This class of interactions is omnipresent in multi-asset systems, where degradation or failure on one asset imposes additional stress on the remaining assets. For instance, degrading gearboxes impose significant stress on shafts and connected bearings.
Degradation models consider multi-asset dependencies through 
\textit{failure-triggered}  or \textit{degradation-triggered} 

interactions. Failure triggered models study cases where a failure of an asset affects either failure or degradation of remaining functioning assets \citep{zequeira2005inspection}.
In degradation-triggered models, degradation interaction occurs continuously throughout asset lifetime. Interacting degradation paths are typically modeled through joint distribution functions of degradation paths 
\citep{wang2004reliability}, stochastic models \citep{bian2014stochastic}, or copula functions \citep{lu2021general}.

Preventive maintenance actions are mostly scheduled periodically in the maintenance literature. Units undergo maintenance in predetermined intervals solely based on failure statistics. Although periodic maintenance is overly conservative and cost-inefficient, it precludes unexpected failure of components, posing a risk to the system reliability \citep{arts2018design}. 
In condition-based maintenance (CBM) policies, maintenance actions are planned for components based on condition monitoring information. 
Most of the existing literature on CBM focuses on single-component systems \citep{yildirim2016sensor}. Although these classes of CBM strategies are adaptable to multi-component systems with independent components \citep{bakir2021integrated}, they have limited applicability in the presence of component dependencies. In line with the literature on degradation modeling of dependent components, literature on CBM scheduling in multi-component systems can be categorized on the premise of stochastic interactions between components, i.e., failure or degradation interactions. The existing literature mainly focused on the failure interactions, i.e., interactions triggered by failure events \citep{sheu2015extended}. 
One drawback of these methodologies is their functionality in practice as the lifetime, and consequently the failure time, of assets are jeopardized by the level of usage, load, or stress they experience \citep{hollander1995dynamic}. The second group of studies aims to model the effect of stochastic interactions at the level of degradation processes. CBM studies on multi-component systems considering stochastic degradation interactions involve complex analytical formulations and offer limited flexibility limited as a result. The degradation-triggered interactions are often considered for small systems \citep{do2015condition} or with limiting assumptions \citep{rasmekomen2016condition}. There are recent studies that model dynamic loading and component-to-component interaction \citep{,yildirim2019leveraging,basciftci2020data}. However, strong assumptions limit their applicability, and they fall short of incorporating both dependencies into O\&M decisions. Fully harnessing degradation models and integrating them into O\&M decisions requires a new generation of models that can embed degradation functions within optimization models and effectively model degradation uncertainties.

A variety of modeling approaches have been used to tackle uncertainty for CBM problems. The literature for them is two-fold. In the first category, there are studies that focus on single-asset systems \citep{alaswad2017review,ding2015maintenance}. In the second category, multi-asset systems have been investigated by approaches such as stochastic programming, chance-constraint approach, Markov decision process, and robust optimization to tackle uncertainty. \cite{feng2015reliability} studies the reliability and CBM of a non-repairable, multi-asset system considering stochastic dependency caused by environmental factors (e.g., temperature). They assume a deterministic system with one major component and several dependent assets. 
In their second work, \cite{feng2016reliability} extended the idea by considering dependence between assets due to randomly arriving shocks. The effect of stochastic dependency on the resulting policy has not been investigated. \cite{mercier2012preventive} and \cite{mercier2014condition} consider a two-component system and describe the dependency through a bi-variate Lévy process. Recently, there has been a growing interest in incorporating sensor-driven prognostic insights into decision optimization models. For handling these issues, some studies incorporate the uncertainty associated with RLDs estimations into the objective function \citep{yildirim2016sensor,bakir2021integrated}. In these studies, a degradation-based cost function is dynamically updated upon the arrival of new sensor observations. Some studies consider the stochastic nature of remaining life estimations with known distributions and utilize the stochastic programming approaches to represent the uncertainties through failure scenarios \citep{basciftci2020data}. 
A major limitation of these approaches relates to poor computational scalability with respect to increasing number of assets and scenarios. 

In this article, we propose a novel methodology for constructing and incorporating the degradation models in a sensor-driven robust optimization framework. The proposed model addressed two-gaps identified in literature. Our approach accounts for the joint consideration of loading impact and multi-asset interaction dependencies in degradation modeling. We offer degradation models that can incorporate both types of dependencies for multi-asset systems. 

The second and more significant gap is integrating these degradation dependencies into a robust decision-optimization model. To the best of our knowledge, we propose the first model that embeds reformulations of decision-dependent degradation signals within optimization models to enable operators to accurately predict, and control asset degradation and failure risks. Proposed methodology does not make any assumptions on distribution of the underlying uncertainties. The aim is to choose the best-immunized O\&M decisions against ``uncertain but bounded" degradation realizations of assets. 
\section{Model Formulation}
\begin{figure*}[htbp!]
    \centering 
    \includegraphics[width=0.9\textwidth]{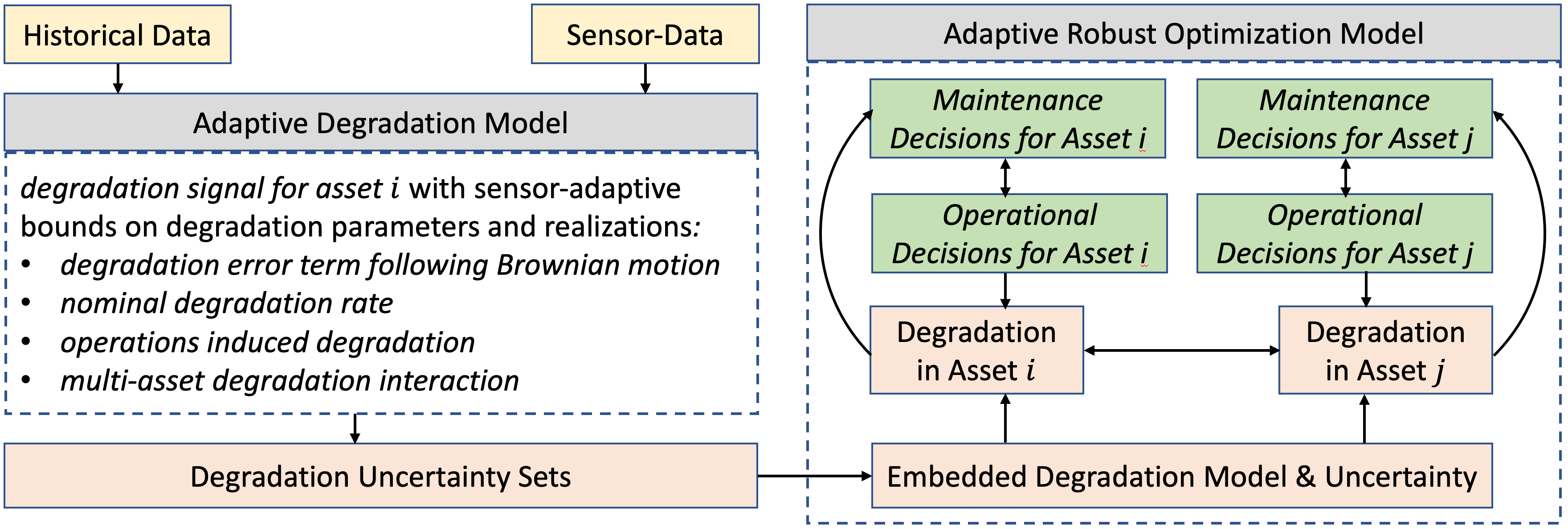}
             \caption{Overview of the Robust Optimization Framework}                      \label{overview}
\end{figure*}

The proposed degradation-driven framework for O\&M builds on a seamless integration across two main components, as shown in Figure 1. The first component of the framework, the degradation model, uses historical and sensor data to develop predictive models to characterize degradation signals in the assets. These degradation signals evolve as a function of three components: (i) inherent degradation rate and error, which defines the progression of degradation in a fixed nominal operating environment, (ii) operations-induced degradation, which models the impact of different operational decisions on asset degradation, and (iii) multi-asset degradation interactions, which capture dependencies across assets in terms of how they degrade and fail. These parameters form the basis for our degradation uncertainty sets. Parameters associated with these uncertainty sets are continuously updated and refined based on streaming sensor data. The second component of the framework, the robust optimization model, embeds sensor-adaptive degradation uncertainty sets and degradation signal formulations within a decision optimization formulation. In Figure \ref{overview}, arrows indicate dependencies across decisions and their outcomes. 

This section formulates different components of the framework and demonstrates how they are integrated. Section \ref{Degradation Modeling and Prognostics} describes the degradation models developed to predict asset failure risks. Section \ref{Reformulating_model} develops a reformulation of the degradation models to embed them within a mathematical programming framework. Using this reformulation, Section \ref{Deterministic model} introduces the deterministic formulation for the proposed O\&M optimization model. Section \ref{Uncertainty Set} discusses different approaches to modeling sensor-driven uncertainty sets, and Section \ref{Revisiting Deterministic Formulation} revisits the deterministic formulation to offer an alternative model. Finally, Section \ref{robust count} presents the robust counterpart of the problem.

\subsection{Multi-Asset Degradation Models}
\label{Degradation Modeling and Prognostics}

Degradation modeling develops stochastic functions to capture the time-varying behavior of degradation signals - i.e., modeling both the current state and probable trajectories of health. Typically, similar assets exhibit the same functional form for their degradation, e.g.,
degradation signals follow an increasing exponential trend over time. However, they exhibit significant differences in the parameters of these functions, such as the rate of degradation. This variability is typically due to nuances in installation, manufacturing, etc. 
When the assets operate in time-varying conditions, the additional stress due to these conditions imposes further variability on the rate of degradation.  
The main objective of the degradation modeling approaches is to provide an accurate model of the degradation trajectory in order to predict when the degradation severity reaches a certain threshold level, which constitutes failure \citep{zhou2021industrial,gebraeel2005residual}. This section elucidates the degradation modeling framework, by introducing factors that impact degradation, developing a degradation model,
and using these formulations to predict asset remaining life.

\subsubsection{Decision-Dependent Degradation Factors:} In our modeling framework, we consider two types of decision-dependent degradation factors that impact the failure risks of a fleet of assets: operations-induced degradation and multi-asset degradation interactions. 

\begin{itemize}
\item \underline{\textit{Operations-induced degradation (OID)}} relates to the dynamic stress on assets as a function of operational decisions. Typically, when assets are operated in more demanding operational regimes (e.g., increased production); they exhibit an increasing pace of degradation. This interaction introduces a significant tradeoff: \textit{Would it be optimal to increase production for operational benefits or ease production to extend equipment lifetime and reduce long-term maintenance costs?}
\item \underline{\textit{Multi-asset degradation interactions (MDI)}} refers to degradation rate interactions across assets. This type of interaction is omnipresent in systems with multiple connected assets, where a degradation in one asset exacerbates the rate of degradation in other connected assets. A second tradeoff arises due to this interaction: \textit{Would it be optimal to replace an asset earlier than planned, if it causes other critical assets to degrade at a significantly higher pace?}
\end{itemize}
These complex tradeoffs resulting from decision-dependent degradation factors are further compounded by their impact on fleet-level O\&M decisions, specifically with regard to the optimal production rate and asset replacement policies. Addressing these intertwined considerations requires a comprehensive modeling framework that can jointly assess all these factors. Detailed information and practical examples of these degradation factors are provided in the online appendix

\subsubsection{Degradation Model:}
We propose a degradation model that inherently captures both OID and MDI degradation factors. In our modeling framework, the degradation signal for asset $i$, $l_i(t)$, is a continuous-time continuous-state stochastic function:
\begin{align} 
{l}_i(t) = {l}_i(0) + \int_{0}^{t} f\Big[s;\underbrace{\vphantom{\int_{0}^{t}}
{\mathcal{D}_i}}_{\substack{\text{Nominal} \\ \text{Rate}}},\underbrace{\vphantom{\int_{0}^{t}}
\zeta_i(\rho_i(s))}_{\substack{\text{OID} \\ \text{Impact}}
},\underbrace{\vphantom{\int_{0}^{t}}
\boldsymbol{\gamma}_i(\boldsymbol{l}(s))}_{\substack{\text{MDI} \\ \text{Impact}}
}\Big]\,ds \  + \epsilon_i(t)
\label{eq:deg_gen}
\end{align}
The proposed degradation signal is composed of four components. The first component is the nominal degradation governed by rate $\mathcal{D}_i$, which denotes the inherent rate of degradation. It is the pace of degradation that occurs independently of our decisions. Second, we introduce OID impact on asset degradation, $\zeta_i(\rho_i(t))$, where $\rho_i(t)$ denotes the operational decision i.e., loading, of asset $i$ taken at time $t$ and $\zeta_i(.)$ is a function that maps operational decision to its associated OID impact. Third, we focus on the MDI interactions: $\boldsymbol{\gamma}_i(\boldsymbol{l}(t)) = \boldsymbol{\gamma}_i({l}_j(t), \; \forall j \in \mathcal{A}_i)$ evolves as a function that evaluates MDI impact by using degradation signals from a set of assets $\mathcal{A}_i$ that interacts with asset $i$. In other words, this interaction explicitly links the rate of degradation in asset $i$ with the severity of degradation in interconnected assets. We would like to note that the second and third components of the model are decision-dependent degradation factors, which enable the decision maker to control their contributions to the degradation rate. Lastly, we model the inherent uncertainty of degradation, which follows a Brownian motion process $\epsilon_i(t) = B_i(t) \sim N(0,\sigma^2_i t)$. 

To ensure that the proposed degradation model can be incorporated into a robust optimization framework, we consider a special class of the degradation model introduced in equation \eqref{eq:deg_gen} using the following four assumptions. 
First, we assume a linear relationship between the degradation rate
of asset $i$ and the degradation 
severity of other interdependent assets $j \in \mathcal{A}_i$
, e.g., an elevation in the degradation signal induces a linearly-increasing behavior in the rate of degradation of other interdependent assets. Second, we assume a linear relationship between the rate of degradation
and loading $\rho_i(t)$ of asset $i$, e.g., raising operational loading triggers an increase in the degradation rate of the component. 
The third assumption is that the degradation rate
is a linear function of its nominal degradation rate, OID and MDI impact factors, i.e., 
$l_i(t)$:=$\mathcal{D}_i +\zeta_i(\rho_i(t))+ \sum_{j\in\mathcal{A}_i}\boldsymbol{\gamma}_{i}(l_j(t))$. As a result, we present the following reformulation:

\begin{equation} 
{l}_i(t) = {l}_i(0) + \int_{0}^{t} \left[\mathcal{D}_i +\zeta_i( \rho_i(s))+ \sum_{j\in\mathcal{A}_i}\gamma_{j,i}(l_j(s))\right]\,ds \  + B_i(t) \label{eq:deg_spec}
\end{equation}
where $\mathcal{D}_i$, $\zeta_i(.)$ and $\gamma_{j,i}(.)$
are stochastic variables. 
Sensor information can be harnessed to update the distributions of these functions in real time using Monte Carlo methods for Bayesian Inference. 
We assume that a failure occurs when an asset's
degradation signal ${l}_i(t)$ exceeds the failure threshold $\Lambda_i$ for the first time. We observe a sequence of degradation signals up to the observation time $t_i^o$, $\{l_1,\dots,l_k\}$, where $l_k$ is observed at time $t_i^o$. Using these observations, we compute the posterior estimates of the degradation parameters, and remaining life distribution of asset $i$ as follows:
{
\begin{align}  \label{eq:ProbRi}
P\{ R^i_{t_i^o} = t| l_1, \dots, l_k \} = P\left\{ t = \min\left\{ s>0 :  {l}_i(s+t_i^o) \geq \Lambda_i \right\}| l_1, \dots, l_k \right\}
\end{align}}

which may either have a closed-form solution or may need to be obtained through numerical methods \citep{gebraeel2005residual}. Next, we will focus on integrating these predictions on remaining life distributions and degradation functions within a decision optimization framework.

\subsection{Reformulating Degradation Model within Mathematical Programming}
A central modeling challenge in our framework is to establish a link across the failure prediction model shown in equation \eqref{eq:ProbRi} and our proposed optimization model. Existing literature focuses on precomputing these failure risks before solving an optimization model, and using the resulting failure risks to generate maintenance cost functions \citep{yildirim2016sensor}, failure scenarios \citep{basciftci2020data}, or Markov decision process formulation \citep{uit2020condition}. These modeling approaches assume a one-directional interaction, i.e., degradation outcomes impact O\&M decisions. In our framework, this interaction is bi-directional, meaning that the degradation also evolves as a decision-dependent function. For any asset $i$, decision dependency is twofold: (i) operational decisions of asset $i$ impacts its rate of degradation, and (ii) operational decisions of asset $j$ impacts the rate of degradation in asset $j$, which in turn impacts the rate of degradation in asset $i$ due to multi-asset degradation interactions. To incorporate these factors into our model, our approach (i) develops a reformulation of the degradation signals in equation \eqref{eq:deg_spec}, and (ii) reconstructs the failure evaluation elucidated in equation \eqref{eq:ProbRi} within the optimization model. 

To integrate the degradation dynamics into an optimization model, we introduce a binary decision variable, $u_i(t)$ that encodes whether unit $i$ is unavailable. The variable should satisfy the following logic (i) if the degradation level $l_i(t)$ exceeds the threshold $\Lambda_i$, the binary variable $u_i(t+1)$ is set to one, otherwise, (ii)  the degradation process evolves as governed by the differences between 
$l_i(t+1)$ and $l_i(t)$, which is derived from equation \eqref{eq:deg_spec}. We impose this relation via linear constraints, discretized over the planning horizon $H$:

\begin{equation}
    l_{i,t}\geq l_{i,0}+  \sum_{\tau = 1}^{t}\bigg( d_{i,\tau} +  \zeta_{i,\tau}. p_{i,\tau}+\sum_{j \in \mathcal{A}_i} \gamma_{j,i,\tau}.\omega_{j,\tau} \bigg)- M  \sum_{\tau = 1}^{t} u_{i,\tau} , \hspace{0.5cm}  \forall i \in \mathcal{A},\hspace{0.05mm}\forall t \in \mathcal{T}  \label{degradation_cons}
\end{equation}

To enable the transition to constraint \eqref{degradation_cons}, we perform the following modifications: 
\begin{inparaenum}[\itshape (i)] \item $l_{i,t}$ within the optimization model represents a time-discretized version of the degradation signal $l_i(t)$ of the degradation model for each asset $i\in \mathcal{A}$, and time $t \in \mathcal{T}$. \item $d_{i,t}$ denotes the summation of the nominal rate of degradation and uncertainty between times $t$ and $t+1$, i.e., $d_{i,t}:=$$\mathcal{D}_i+(\epsilon_i(t+1)-\epsilon_i(t))$. Using the incremental property of Brownian processes, we can conclude that the error increment $(\epsilon_i(t+1)-\epsilon_i(t))$ follows a normal distribution. \item Loading impact $\zeta_i(\rho_i(t))$ is a function of operational loading. In our formulation, we introduce time-discretized $p_{i,t}$ and $\zeta_{i,t}$ to denote the operation and its contribution to degradation, respectively. The operational decision in our framework is the production level.\item Similarly, we define time-discretized $\omega_{j,t}$ and $\gamma_{j,i,t}$ as the degradation level of asset $j\in\mathcal{A}_i$ at time $t$ and its corresponding impact on asset $i$, respectively. Finally, \item the availability status of asset $i$ is represented by the binary variable $u_{i,t}$. 
\end{inparaenum}

We note that the availability of the asset makes a significant impact on the degradation level. In our framework, unavailability may occur due to two reasons, maintenance or failure. Maintenance may be scheduled for asset $i$ at time $t$, which sets the unavailability due to the maintenance variable to $u^m_{i,t}$ = 1; or a failure may occur, which sets the unavailability due to the failure variable to $u^f_{i,t}$ = 1. Evidently, $u_{i,t} = u^m_{i,t}+u^f_{i,t}$ for each $i\in\mathcal{A}$ and $t\in\mathcal{T}$. Whenever an asset becomes unavailable,
degradation level is reset to zero, meaning the asset is to become brand-new.
This is achieved via (i) the use of big-M formulation that relaxes the lower bound of the degradation signal to zero, hence (ii) objective function incentivizing the degradation level to go down to its minimum level.

\label{Reformulating_model}

\subsection{Deterministic Formulation Considering Single Maintenance}
The degradation model reformulation in Section \ref{Reformulating_model} is integrated within a large-scale O\&M model for a fleet of assets. Decision variables for the proposed model is introduced in Table \ref{tab:vars1}. Our maintenance decisions identify both preventive \textit{(before-failure)} and corrective \textit{(after-failure)} maintenance actions for each asset. Both maintenance actions lead to asset downtime, which interrupts production throughput. Corrective maintenance cost $C^c$ is typically greater than preventive maintenance cost $C^p$. A failed asset remains out of service until corrective maintenance is completed. 
. 

Operational decisions include production variables that collectively determine the fleet-level throughput. The parameter $G_{i,t}$ denote the unit production cost of asset $i$ at time $t$ and per unit of unsatisfied demand is penalized by $C^u$. Building upon the degradation model incorporated through constraints \eqref{degradation_cons}, we enable bidirectional interactions across O\&M decisions. On one hand, unavailability of assets due to failure or maintenance outages impacts the throughput.

On the other hand, operational decisions impact the degradation rates, failure risks, and maintenance schedules.

\begin{table}
\centering
\label{table: vardescription1}
  \caption{Decision Variables for the Deterministic Formulation} 
  \renewcommand{\arraystretch}{0.65}
  \begin{tabular}{ll}
  \hline
  Variable & Description \\
		\hline
		$m^p_{i,t} \in \{0,1\}$ \quad & \text{$m^p_{i,t}$=1, if a preventive maintenance for asset $i$ starts at time $t$} \\
		$m^c_{i,t} \in \{0,1\}$ \quad & \text{$m^c_{i,t}$=1, if a corrective maintenance for asset $i$ starts at time $t$} \\
		$u^m_{i,t} \in \{0,1\}$ \quad & \text{$u^m_{i,t}$=1, if asset $i$ is unavailable
  due to maintenance at time $t$} \\
		$u^f_{i,t} \in \{0,1\}$ \quad & \text{$u^f_{i,t}$=1 if asset $i$ is unavailable
  due to failure at time $t$} \\
		$p_{i,t} \in \mathbb{R}_{+}$ \quad & \text{Production level of asset $i$ at time $t$}\\
		$\psi_{t} \in \mathbb{R}_{+}$ \quad & \text{Unsatisfied demand at time $t$}\\
		\hline 
		\end{tabular}
  \label{tab:vars1}
	\end{table}
\begin{inparaenum}[\itshape (i)] 
\underline{\textit{Objective}} of the decision maker is to optimize fleet-level O\&M decisions that \item leverage sensor-driven degradation models to conduct maintenance when needed and reduce the risks of unexpected failures, \item minimize the impact of unavailability on operations, and \item finetune fleet-level stress/loading on the assets to control the rate of degradation. These factors are reflected within the objective function, which minimizes the total O\&M cost.
\end{inparaenum}
\begin{align}
&Min \quad \displaystyle\sum_{i \in \mathcal{A}}\displaystyle \sum_{t \in \mathcal{T}} (C^{p} \cdot m^p_{i,t}  +C^{c}\cdot m^c_{i,t}) +\sum_{i \in \mathcal{A}}C^{c}\cdot u^f_{i,H} + \sum_{i \in \mathcal{A}}\displaystyle \sum_{t \in \mathcal{T}} G_{i,t} \cdot p_{i,t}+\sum_{t \in \mathcal{T}}C^{u}\cdot \psi_{t}  \hspace{1cm}\label{objective function}
\end{align}
The first two terms represent preventive and corrective maintenance costs incurred over the planning horizon, respectively. The third term charges corrective maintenance costs unless the failed assets are maintained before the end of the planning horizon $H$. The last two terms represent the cost of production and unsatisfied demand over the planning horizon, respectively.

\underline{\textit{Degradation constraints}} are enforced along with constraints \eqref{degradation_cons} developed in Section 3.2, to model the degradation behavior. Constraints \eqref{component_0.1} - \eqref{component_2} account for updating degradation levels of interacting asset pairs. As mentioned earlier, we assume that asset $j$'s contribution to asset $i$'s degradation rate is a linear function of asset $j$'s degradation level. Specifically, the higher the degradation in asset $j$, the more it accelerates the degradation rate in asset $i$. There are two exceptions to this rule: \begin{inparaenum}[\itshape (1)]
  \item when asset $j$ fails, we assume that the impact on asset $i$ reaches its maximum value (i.e., proportional to the failure threshold $\Lambda_i$), and \item when asset $j$ is undergoing planned preventive maintenance, we assume there is no interaction between asset $j$ and $i$.
\end{inparaenum}  These exceptions are the main reason for introducing the additional variable $\omega_{i,t}$, which is identical to the degradation level of asset $i$, $l_{i,t}$ except for the aforementioned exceptions. 
\vspace{-2mm}
\begin{flalign}
&u_{i,t} = u^m_{i,t}+u^f_{i,t}, &\forall i \in \mathcal{A},\hspace{1mm} \forall t \in \mathcal{T}\label{component_0.1} \\
&l_{i,t} \leq \Lambda_i, &\forall i \in \mathcal{A},\hspace{1mm} \forall t \in \mathcal{T}\label{component_0.2}\\
&  \omega_{i,t} \geq l_{i,t}, 
& \forall i \in \mathcal{A},\hspace{1mm}\forall t \in \mathcal{T} \label{component_2.0}\\
& \omega_{i,t} \geq \Lambda_i \cdot u^f_{i,t},   &\forall i \in \mathcal{A},\hspace{1mm}\forall t \in \mathcal{T}
\label{component_2}
\end{flalign}
\vspace{-1mm}
Constraints \eqref{component_0.1} ensure that unavailability occurs either due to preventive maintenance or failure, and constraints \eqref{component_0.2} ensure that the degradation signal is bounded by the failure threshold. Constraints \eqref{component_2.0} and \eqref{component_2} define the lower bounds for the interaction term $\omega_{i,t}$: degradation amplitude always constitutes a lower bound, but an additional lower bound of $\Lambda_i$ is enforced when the asset is unavailable due to failure.

. 
\vspace{-1mm}
\underline{\textit{Maintenance Constraints}} collectively model preventive and corrective maintenance actions, and their impacts on asset availability and the use of maintenance crew resources. 
\begin{flalign}
& u_{i,t}^m  = \sum_{\tau = max\{0, t - Y^p +1\}}^{t}m_{i,\tau}^p + \sum_{\tau = max\{0, t - Y^c +1\}}^{t}m_{i,\tau}^c && \forall i \in \mathcal{A},\hspace{1mm}\forall t \in \mathcal{T}\label{11} \\
&1-u^f_{i,t} \geq m^p_{i,t} && \forall i \in \mathcal{A},\hspace{1mm}\forall t \in \mathcal{T}\label{12} \\
&   u^f_{i,t-1}- u^f_{i,t} \leq  m^c_{i,t} && \forall i \in \mathcal{A},\hspace{1mm}\forall t \in \mathcal{T}\label{13} \\
&\displaystyle\sum_{i\in\mathcal{A}}  u^m_{i,t} \leq Q    &&\hspace{1.2cm} \forall t \in \mathcal{T} \label{14}
\end{flalign}

\vspace{-3mm}
Constraint \eqref{11} ensures the asset remains unavailable during maintenance, where $Y^p$ and $Y^c$ denote preventive and corrective maintenance durations, respectively. If a preventive maintenance starts between $t-Y^p+1$ and $t$, we conclude that there must be an ongoing maintenance. Similar logic applies for the corrective maintenance actions. Constraint \eqref{12} prohibits conducting preventive maintenance when an asset is failed. Constraint \eqref{13} ensures that a corrective maintenance is initiated when $u^f_{i,t}$ variable switches from 1 to 0. In Constraint \eqref{14}, we restrict the number of ongoing simultaneous maintenance to predetermined maintenance crew capacity Q.
\\\underline{\textit{Operational Constraints}} model the maintenance impact on asset production, and fleet throughput.
\vspace{-3mm}
\begin{flalign}
&\displaystyle  p_{i,t} \leq K_i \cdot( 1- u_{i,t} )   &\forall i \in \mathcal{A},\forall t \in \mathcal{T}\label{15}\\
&\sum_{i \in \mathcal{A}} p_{i,t} + \psi_{t} \geq S_{t}& \forall t \in \mathcal{T}\label{17}
\end{flalign}
Constraint \eqref{15} establishes the dependence between production throughput and asset downtime. We enforce production suspension if an asset is out of service due to maintenance or failure. Otherwise, the asset can produce up to its production capacity $K_i$ when the asset is available to produce.

Constraint \eqref{17} counts the unsatisfied demand as the difference between the demand $S_t$ and
fleet-level production throughput for each time $t$. Finally, constraint \eqref{18} constitutes the sign and binary restrictions.
\begin{flalign}
&l_{i,t}, \omega_{i,t},p_{i,t},\psi_t \in {\mathbb{R^{+}}},\hspace{1cm}
m^c_{i,t},m^p_{i,t}, u_{i,t}, u^f_{i,t},u^m_{i,t}\in \{ 0,1\},\hspace{0.5cm} \forall i \in \mathcal{A},\; \forall t \in \mathcal{T} \label{18}
\end{flalign}
\label{Deterministic model}

\vspace{-6mm}
\label{Nested uncertainty set section}
\subsection{Deterministic Formulation Considering Multiple Maintenances}

\begin{table}
\centering\label{table: vardescription2}
  \caption{Additional Decision Variables} 
  \renewcommand{\arraystretch}{0.65}
  \begin{tabular}{ll}
  \hline
    Variable & Description \\
    \hline
        $z_{i,t,k} \in \{0,1\}$  \quad &  \text{$z_{i,t,k}$ =1, if time $t$ is not between end time of ${(k-1)}^{th}$} and start time of $k^{th}$\\
		& \text{maintenance of asset $i$} \\
		$v_{i,t,k} \in \{0,1\}$ \quad & \text{$v_{i,t,k}$ =1, if $k^{th}$ maintenance for asset $i$ starts at time $t$}\\
		$v_{i,k}^{0} \in \{0,1\}$ \quad & \text{$v_{i,k}^{0}$ =1, if $k^{th}$ maintenance for asset $i$ is not scheduled}\\ 
		$l^{'}_{i,t,k} \in \mathbb{R}_{+}$ \quad & \text{Total degradation level for asset $i$ at time $t$ before $k^{th}$ maintenance} \\
		$\omega_{j,i,t,k}^{'} \in \mathbb{R}_{+}$ \quad & \text{Degradation impact of asset $j$ on asset $i$ at time $t$ before $k^{th}$ maintenance}\\ & \text{of asset $i$}\\
		$p_{i,t,k}^{'} \in \mathbb{R}^{n}_{+}$ \quad & \text{Operational loading on asset $i$ at time $t$ before $k^{th}$ maintenance}\\
		\hline 
		\end{tabular}\label{tab:vars2}
\end{table}
This section presents a reformulation of the deterministic model discussed in Section \ref{Deterministic model} to accommodate multiple maintenance actions within planning horizon.
We use a representation of maintenance cycles $\mathcal{K}$, where $k=1$ refers to the first maintenance cycle, $k=2$ refers to the second maintenance cycle, and so forth. In this approach, each maintenance event activates a new maintenance cycle $k$, deactivates the previous maintenance cycle $k-1$, and resets the total accumulated degradation level to zero. To keep track of the maintenance cycles, a new set of decision variables is introduced (see Table \ref{tab:vars2}). {\ A sample maintenance schedule is provided in online appendix to illustrate the interaction across decision variables $v$, $z$ and $u^m$.}

The new formulation uses with a cumulative function that includes degradation increments from time $0$ to time $t$ and considers any maintenance actions that reset the degradation cycle. The following constraints define cumulative degradation at time $t$:   
\begin{align}
& l^{'}_{i,t,k} \geq l^{'}_{i,0,k} + \sum_{\tau = 1}^{t}\bigg(d_{i,\tau}(1-z_{i,\tau,k}-u_{i,\tau}^{f}) +  \zeta_{i,\tau}. p_{i,\tau,k}^{'}+\sum_{j \in \mathcal{A}_i} \gamma_{j,i,\tau}.\omega_{j,i,\tau,k}^{'}\bigg)\nonumber \\ & \hspace{7.5cm}\forall i \in \mathcal{A},\forall t \in \mathcal{T},\forall k \in \mathcal{K} \label{degradation sum} \\
& l^{'}_{i,t,k} \leq \Lambda_i, \hspace{6.5cm} \forall i \in \mathcal{A},\forall t \in \mathcal{T},\forall k \in \mathcal{K}\label{threshold2}
\end{align}

More specifically, constraint \eqref{degradation sum} is another representation of constraint \eqref{degradation_cons} in which the degradation level $l^{'}_{i,t,k}$ is the accumulative degradation that occurs during maintenance cycle $k$. The inherent degradation multiplier $(1-z_{i,t,k}-u_{i,t}^{f})$ is equal to 1 only if time $t$ falls within maintenance cycle $k$ (i.e. $z_{i,t,k} = 0$) and the asset is operational (i.e. $u_{i,t}^{f}=0$). Decision-dependent multipliers $p_{i,t,k}^{'}$ and $\omega_{j,i,\tau,k}^{'}$ are for the OID and MDI loading conditions, respectively. These multipliers are constrained to be nonzero only when time $t$ is within maintenance cycle $k$, and are modeled using the following set of constraints:

\underline{\textit{Constraints for OID impact}} are restated by linking the $p_{i,t,k}^{'}$, the OID impact during maintenance cycle $k$, and production variable $p_{i,t}$ as follows: 
\vspace{-3mm}
\begin{flalign}
& p_{i,t,k}^{'} \geq p_{i,t} - M\cdot z_{i,t,k}&\forall i \in \mathcal{A},\hspace{1mm}\forall t \in \mathcal{T},\hspace{1mm}\forall k \in \mathcal{K}\label{loading}
\end{flalign}
\vspace{-1mm}
{Constraint \eqref{loading} ensures that if $t$ falls within maintenance cycle $k$ (i.e., $z_{i,t,k}=0$), the OID impact multiplier $p_{i,t,k}^{'}$ equals the production level $p_{i,t}$. If $t$ is not within cycle $k$ (i.e., $z_{i,t,k}=1$), $p_{i,t,k}^{'}$ can take any positive value. The constant $M$ is a big-M coefficient.}

\underline{\textit{Constraints for MDI impact}} follows similar logic to reformulate by creating the link between the $\omega_{i,t}$ and $\omega_{j,i,t,k}^{'}$ variables representing the degradation interactions between assets.
\vspace{-6mm}

\begin{flalign}
& \omega_{i,t} \geq l^{'}_{i,t,k} - M\cdot z_{i,t,k} & \forall i \in \mathcal{A},\hspace{1mm}\forall t \in \mathcal{T},\hspace{1mm}\forall k \in \mathcal{K}\label{component_1_updated}\\
& \omega_{j,i,t,k}^{'} \geq \omega_{j,t}- M\cdot z_{i,t,k}- M\cdot  u_{i,t}^{f}&\forall i \in \mathcal{A},\hspace{1mm}\forall j \in \mathcal{A}_i,\hspace{1mm}\forall t \in \mathcal{T},\hspace{1mm}\forall k \in \mathcal{K}\label{component_3_updated}
\end{flalign}
Constraints \eqref{component_1_updated} guarantee that $\omega_{i,t}$ takes the same value with $l_{i,t,k}$ if time $t$ is in the $k^{th}$ maintenance cycle. In constraints \eqref{component_3_updated}, we ensure that asset $j$ does not make any contribution to the degradation increment of its pair $i$ while it is unavailable due to maintenance or failure.

\underline{\textit{Constraints for linking  maintenance, cycles, and unavailability}} are included to ensure that maintenance cycles are coordinated with the unavailability and maintenance decisions of the assets.
\begin{flalign}
&\sum_{k \in \mathcal{K}} v_{i,t,k}  \leq 1 &  \forall i \in \mathcal{A},\hspace{1mm}\forall t \in \mathcal{T}\label{5}\\
& v_{i,k-1}^{0} \leq v_{i,k}^{0}  &  \forall i \in \mathcal{A},\hspace{1mm}\forall k \in \mathcal{K}\label{7}\\
&\sum_{t \in \mathcal{T}} v_{i,t,k} +v_{i,k}^{0} = 1 &  \forall i \in \mathcal{A},\hspace{1mm}\forall k \in \mathcal{K}\label{8}\\
&v_{i,t,k}  \leq m_{i,t}^{p} &  \forall i \in \mathcal{A},\hspace{1mm}\forall t \in \mathcal{T},\hspace{1mm}\forall k \in \mathcal{K}\label{6}\\
& z_{i,t,1}\leq \sum_{\tau=0}^{t} v_{i,\tau,1}+ u_{i,t}^{m}&  \forall i \in \mathcal{A},\hspace{1mm}\forall t \in \mathcal{T}\label{9}\\
& z_{i,t,k}\leq (1 -\sum_{\tau=0}^{t} v_{i,\tau,k-1}) +\sum_{\tau=0}^{t} v_{i,\tau,k} + u_{i,t}^{m}&  \forall i \in \mathcal{A},\hspace{1mm}\forall t \in \mathcal{T},\hspace{1mm}\forall k \in \mathcal{K}\label{10}
\end{flalign}
Constraints \eqref{5} ensure that each asset can have at most one maintenance cycle in each time period.  Constraints \eqref{7} ensure that maintenance cannot be scheduled for cycle $k$ unless cycle $(k-1)$ has been scheduled, enforcing a logical sequence across maintenance cycles. Constraint set \eqref{8} guarantees that if a starting time of the $k^{th}$ maintenance for asset $i$ is chosen, i.e., $v_{i,t,k}=1$, then $k^{th}$ maintenance must be scheduled, i.e., $v^0_{i,k}=0$. Constraints \eqref{6} mandate that the maintenance variable $m_{i,t}^p$ must be equal to 1 if maintenance of asset $i$ during cycle $k$ starts at time $t$. In constraints \eqref{9} and \eqref{10}, we ensure that $z_{i,t,k}$ must be zero when time $t$ is in between $k^{th}$ and $(k-1)^{th}$ maintenance.

\subsection{Uncertainty Set}
In this section, we introduce the proposed uncertainty set, which 
revolves around three parameters to represent the underlying cause of degradation uncertainty, namely:\begin{inparaenum}[\itshape (i)]\item
	inherent degradation increment  $\boldsymbol{d}$
	, \item OID impact $\boldsymbol{\zeta}$ that accounts for the operations-induced degradation, and \item MDI impact $\boldsymbol{\gamma}$ that captures the impact of multi-asset interactions.  \end{inparaenum}. 

A unique aspect of our uncertainty sets is that it seamlessly adapts to the parameter distributions predicted though the sensor data.

\begin{figure}[h!]
	\centering
	\includegraphics[width=0.65\textwidth]{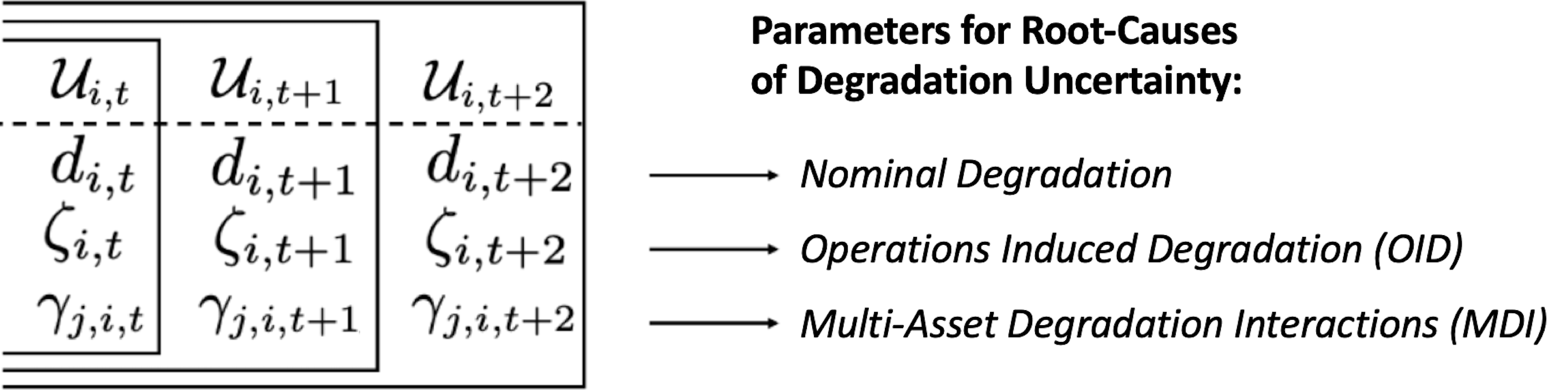}
	\caption{The proposed nested uncertainty set formulation for degradation parameters}
	\label{fig:uncertaintyset}
\end{figure}
\label{Uncertainty Set}

The proposed nested uncertainty jointly models the uncertainties within $t \cdot \left(2+|\mathcal{A}_i| \right)$ parameters, which enables the robust optimization model to consider credible, non-trivial, and realistic realizations of degradation uncertainty. Figure \ref{fig:uncertaintyset} provides an overview of the proposed nested uncertainty set formulation that uses budgeted uncertainty set \citep{bertsimas2006robust}, and the corresponding formulation is given below:
\begin{equation}\mathcal{U}_{i,t}=\left\{
	\begin{aligned}&\sum_{\tau = 1}^{t} \frac{d_{i,\tau}-\bar{d}_{i,\tau}}{\hat{d}_{i,\tau}}  + \sum_{\tau = 1}^{t}\frac{\zeta_{i,\tau}-\bar{\zeta}_{i,\tau}}{\hat{\zeta}_{i,\tau}} \\&\qquad \qquad \qquad \qquad \qquad   + \sum_{\tau = 1}^{t}\sum_{j \in \mathcal{A}_i} \frac{\gamma_{j,i,\tau}-\bar{\gamma}_{j,i,\tau}}{\hat{\gamma}_{j,i,\tau}} \leq  \Delta_{i,t},\\
		&d_{i,\tau} \in [{\bar{d}}_{i,\tau}-\hat{d}_{i,\tau},{\bar{d}}_{i,\tau}+\hat{d}_{i,\tau}], \hspace{1.5cm}\forall \tau \in [1,t]\\
		&\zeta_{i,\tau} \in [\bar{\zeta}_{i,\tau} -\mathcal{\hat{\zeta}}_{i,\tau},{\mathcal{\bar{\zeta}}}_{i,\tau}+\mathcal{\hat{\zeta}}_{i,\tau}], \hspace{1.65cm}\forall \tau \in [1,t]\\
		&\gamma_{j,i,\tau}\in [\bar{\gamma}_{j,i,\tau} - \mathcal{\hat{\gamma}}_{j,i,\tau},{\mathcal{\bar{\gamma}}}_{j,i,\tau}+\mathcal{\hat{\gamma}}_{j,i,\tau}],  \hspace{0.55cm}\forall j \in \mathcal{A}_i, \forall \tau \in [1,t]
	\end{aligned}\right\}\hspace{0.01cm} \forall i \in \mathcal{A}, \hspace{1mm}\forall t \in \mathcal{T}\label{uncertainty_set22}
\end{equation}

In each uncertainty set {$\mathcal{U}_{i,t}$}, we restrict the cumulative deviations from the nominal values of all three stochastic parameters $d_{i,t}$, $\zeta_{i,t}$, $\gamma_{j,i,t}$ up to time $t$. Note that we do not stipulate the agreement of realizations of a parameter in different subsets. In other words, the realization of the stochastic parameters in different subsets can be different, as long as the total deviation of all parameters up until time $t$ is within the specified budget $\Delta_{i,t}$. To illustrate, let $d_{i,t}^{t^1} \in \mathcal{U}_{i,t^1}$ and $d_{i,t}^{t^2} \in \mathcal{U}_{i,t^2}$ denote the realizations of $d_{i,t} $ from $\mathcal{U}_{i,t^1}$, $\mathcal{U}_{i,t^2}$, respectively. We allow realizations, whereby $d_{i,t}^{t^1} \neq d_{i,t}^{t^2}$  for ${t^1} \neq {t^2}$. 

Restricting the summation of deviations over time periods and allowing distinct realizations for the same parameters enable us to capture a wider variety of worst-case realizations within the uncertainty set and yield more robust solutions through the optimization model.

\label{Revisiting Deterministic Formulation}
\subsection{Robust Counterpart} 
\label{robust count}
The proposed extended deterministic formulation does not account for the uncertainty in degradation parameters. However, this assumption may not always hold in practice, as degradation is often subject to variability and uncertainty. Neglecting this uncertainty can lead to unexpected {and costly} failures, which would have a dramatic impact on operational outcomes. To address this issue,
we propose a robust counterpart formulation that quantifies and models this degradation uncertainty in order to ensure that the solutions are immune to any degradation realization within the specified uncertainty set.

To this end, we present the robust counterpart of Constraint \eqref{degradation sum} as follows:
\begin{flalign}
&l^{'}_{i,t,k}-l^{'}_{i,0,k}\geq \nonumber\\&\max_{\boldsymbol{d},\boldsymbol{\zeta},\boldsymbol{\gamma} \in \mathcal{U}_{i,t}}\hspace{2mm} \left( \sum_{\tau = 1}^{t}d_{i,\tau}(1-z_{i,\tau,k}-u_{i,\tau}^{f}) + \sum_{\tau = 1}^{t}\zeta_{i,\tau}. p_{i,\tau,k}^{'}+ \sum_{\tau = 1}^{t}\sum_{j \in \mathcal{A}_i} \gamma_{j,i,\tau}.\omega_{j,i,\tau,k}^{'}\right) \nonumber\\&\hspace{7.9cm} \forall i \in \mathcal{A},\hspace{2mm} \forall t \in \mathcal{T},\hspace{2mm}\forall k \in \mathcal{K}
\label{robust_degradation}
\end{flalign} 

Here, $\mathcal{U}_{i,t}$ denotes the nested degradation uncertainty set for asset $i$ at time $t$, and $\boldsymbol{d}$, $\boldsymbol{\zeta}$, and $\boldsymbol{\gamma}$ represent the degradation increments, OID, and MDI factors, respectively. In the robust counterpart formulation, we introduce an inner optimization problem that maximizes the degradation 
under the uncertainty set for each asset $i$ and time $t$. 
This inner problem serves as a worst-case scenario analysis that allows us to ensure that the solutions are immune to any degradation realization within the uncertainty set. Next, we need to incorporate the worst-case scenario degradation 
into our main optimization problem. This can be achieved by dualizing constraint set \eqref{robust_degradation} and reformulating it as a minimization problem, as shown in Lemma 1. With the uncertainty set \eqref{uncertainty_set22} expressed as a set of constraints, we can then formulate the inner optimization problem as follows:

\begin{lemma} Constraint set \eqref{robust_degradation} can be dualized and reformulated as a minimization problem, which leads to the final reformulation of the proposed robust model as follows: $Min \{ \;(\ref{objective function}) \; |  \; s.t. \;  (\ref{component_2}-\ref{18}),(\ref{threshold2}-\ref{10}) \nonumber, \mathcal{L} \geq \boldsymbol{\pi} \mathbf{u} ,\boldsymbol{\pi}\mathbf{d} \geq \mathbf{m} \}$, where $\mathcal{L}$ denotes left hand side of constraint set \eqref{robust_degradation}, $\boldsymbol{\pi}$ represents corresponding dual variables for uncertainty set constraints, $\mathbf{d}$ and $\mathbf{u}$ represent uncertainty set parameters and $\mathbf{m}$ denotes maintenance and degradation decision variables.
\end{lemma}

Proof of Lemma 1 is provided in the online appendix. The resulting formulation is a robust optimization model that ensures the solutions are immune to worst-case degradation realizations.
\section{Computational Experiments}

In this section, we demonstrate the performance of the proposed O\&M model through a comprehensive case study. The study builds on vibration based degradation data and failure instances acquired from a rotating machinery application \cite{gebraeel2005residual} to determine the degradation parameters, and uses these parameters to emulate degradation processes subjected to a range of degradation configurations. We also demonstrate the computational performance of the proposed model for different problem instances using an acceleration method.

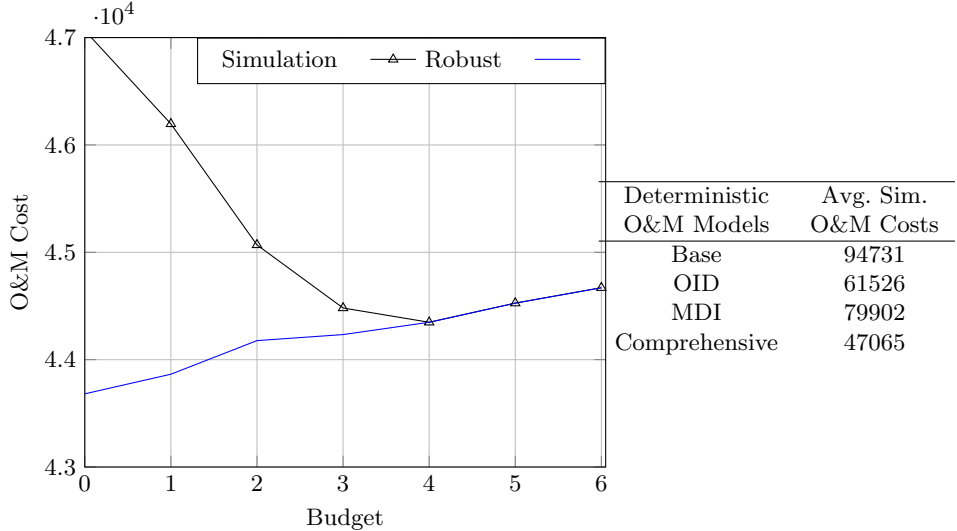
\begin{figure}
\centering
\begin{minipage}{0.6\textwidth}\small
   \centering
    \begin{tikzpicture}
    \begin{axis}
    [name=plot, xlabel={Budget},ylabel={O\&M Cost},
    ymin=43000,ymax=47000, xmin =0, xmax = 6.05,grid=both]
    \addplot[color = black, mark = triangle] table{./Both_Simulation.txt};\label{Simulation}
    \addplot[color = blue] table{./Both_Robust.txt};\label{Robust}
    \end{axis}
    \node[anchor=north  east,fill=white,draw=black] (legend) at ($(plot.north east)$) {\begin{tabular}{ l l l l}
        Simulation & \ref{Simulation}  
        Robust & \ref{Robust}\\
    \end{tabular} };
    \end{tikzpicture}
\end{minipage}%
\begin{minipage}{0.3\textwidth} \small
    \begin{tabular}{cc} \hline
  Deterministic & Avg. Sim. \\ 
  O\&M Models & O\&M Costs \\ \hline
  Base & 94731 \\
  OID  & 61526 \\
  MDI & 79902 \\
  Comprehensive & 47065 \\
  \end{tabular}
\end{minipage}
\caption{Simulation vs Robust Solutions for Comprehensive Model \& Deterministic Models Simulation Results}
\label{fig:synthetic_data_experiments_1}
\end{figure}

\begin{figure}
\centering
\begin{minipage}{0.65\textwidth}
    \centering
    \begin{tikzpicture}
    \begin{axis}
    [name=plot, xlabel={Budget},ylabel={Incurred Penalty Costs},
    ymin=4800,ymax=6100, xmin=0, xmax=6.05,grid=both]
    \addplot[color=black, mark=triangle] table{./Both_Penalty.txt};
    \label{both1}
    \end{axis}
    \node[anchor=north east,fill=white,draw=black] (legend) at ($(plot.north east)$) {\begin{tabular}{ l l l l}
        Simulation & \ref{both1}\\
    \end{tabular}};
    \end{tikzpicture}
\end{minipage}%
\begin{minipage}{0.35\textwidth} \small
    \begin{tabular}{cc}
    \hline
    Deterministic & Avg. Penalty \\ 
  O\&M Models & Cost \\ \hline
  Base & 29882 \\
  OID  & 11270 \\
  MDI  & 28048 \\
  Comprehensive  & 5097 \\
    \end{tabular}
\end{minipage}
\caption{Simulation Outcomes for Penalty Cost \& Deterministic Models Penalty Simulation Results}
\label{fig:synthetic_data_experiments_2}
\end{figure}
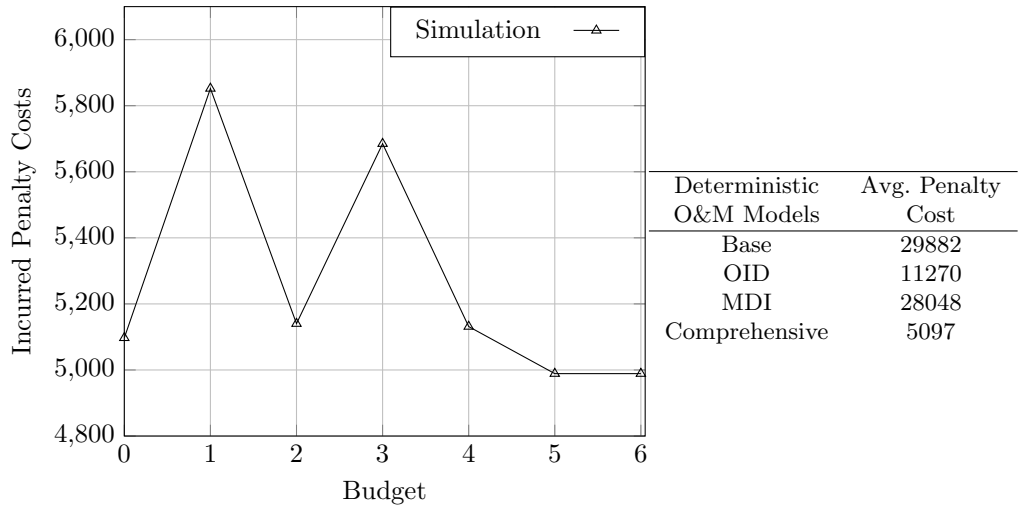

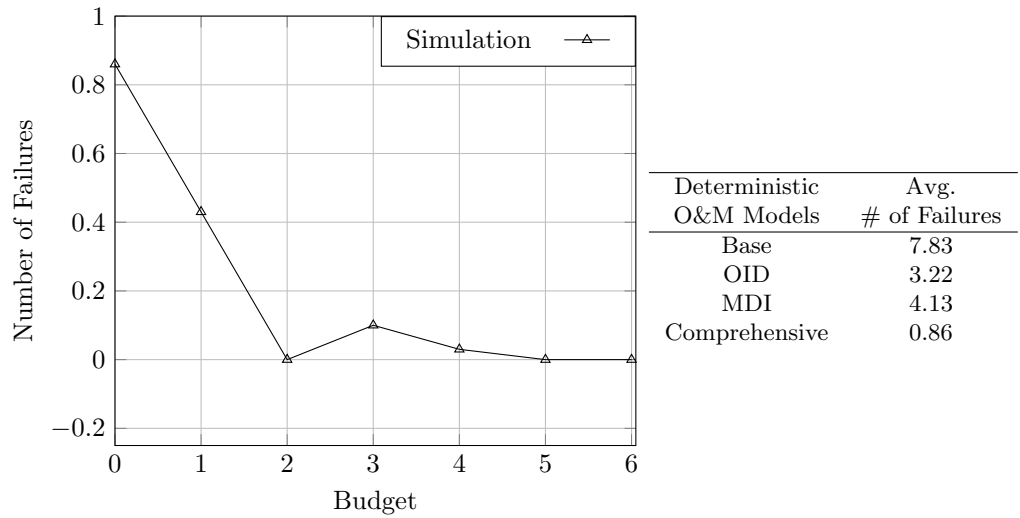
\begin{figure}
\centering
\begin{minipage}{0.65\textwidth}
   \centering
    \begin{tikzpicture}
    \begin{axis}
    [name=plot, xlabel={Budget},ylabel={Number of Failures},
    ymin=-0.25,ymax=1, xmin =0, xmax = 6.05,grid=both]
    \addplot[color = black, mark = triangle] table{./Both_Failure.txt};\label{both}
    \end{axis}
    \node[anchor=north  east,fill=white,draw=black] (legend) at ($(plot.north east)$) {\begin{tabular}{ l l l l}
        Simulation& \ref{both}\\
    \end{tabular} };
    \end{tikzpicture}%
\end{minipage}%
\begin{minipage}{0.35\textwidth} \small
  \begin{tabular}{cc} \hline
    Deterministic & Avg. \\ 
  O\&M Models & \# of Failures \\ \hline
  Base & 7.83 \\
  OID  & 3.22 \\
  MDI  & 4.13 \\
  Comprehensive  & 0.86 \\
  \end{tabular}
\end{minipage}
\caption{Simulation Outcomes for Failure instances \& Deterministic Models Failure Simulation Results}
\label{fig:synthetic_data_experiments_3}
\end{figure}

The case studies are conducted on an experimental framework that brings together two main modules: optimization and degradation simulation. In the optimization stage, we obtain O\&M decisions by solving the optimization model corresponding to each O\&M policy. Then, we fix these decisions and simulate the impact of these decisions under 100 different degradation scenarios. The simulation module implements the decisions from the optimization model as a function of the realized asset failure instances. In the event of failure, the simulation module either chooses to replace the asset or wait. If failure does not occur, the simulation module enforces optimization decisions and evaluates its results. Figure \ref{fig:synthetic_data_experiments_1}, \ref{fig:synthetic_data_experiments_2}, \ref{fig:synthetic_data_experiments_3} shows objective function values of the robust optimization problem, and simulation averages for O\&M costs, penalty costs, and the number of failure instances as a function of different budget values.

While evaluating the performance of the proposed robust optimization model, we consider the impact of different uncertainty sets. The optimization module produces immunized solutions against any realizations within the uncertainty set. The size of the uncertainty set depends on the budget parameter $\boldsymbol{\Delta}$. On one hand, enlarging the uncertainty set provides more protection against uncertainty and mitigates the disruptive impacts of unexpected failures. On the other hand, it elevates the cost of O\&M. We perform experiments with different budget parameters to showcase this trade-off. In addition to budget parameter experiments, we also consider the impact of robust vs deterministic formulations to solve the proposed model. Our case study considers four different policies for managing O\&M decisions:
\begin{itemize}
    \item \textit{Benchmark Policy 1 - Base O\&M Model:} This model considers nominal degradation and ignores the impact of operational loading (OID) and multi-asset degradation interactions (MDI). This is represented by the assumptions $\{ \zeta_{i,t} = 0, \hspace{3mm} \forall i \in \mathcal{A}, t \in \mathcal{T}\}$ and $\{ \omega^{'}_{j,i,t,k} = 0, \hspace{3mm} \forall i \in \mathcal{A},j \in \mathcal{A}_i, t \in \mathcal{T}, k \in \mathcal{K} \}$.
    \vspace{0.1cm}
    \item \textit{Benchmark Policy 2 - O\&M Model with OID:} This model considers both
    nominal degradation and operational loading (OID), but ignores multi-asset degradation interactions (MID), assuming the following holds true $\{ \omega^{'}_{j,i,t,k} = 0, \quad \forall i \in \mathcal{A},j \in \mathcal{A}_i, t \in \mathcal{T}, k \in \mathcal{K} \}$.
    \item \textit{Benchmark Policy 3 - O\&M Model with MDI:} This model considers both 
    nominal degradation and multi-asset degradation interactions (MID), while ignoring operations induced degradation impact (OID). This is represented by the assumption
    $\{ \zeta_{i,t} = 0, \quad \forall i \in \mathcal{A}, t \in \mathcal{T}\}$.
    \item \textit{Proposed Policy - Comprehensive O\&M Model:} This model incorporates all three factors: nominal degradation, operations-induced loading/stress (OID), and multi-asset degradation interactions (MDI).
   With no assumptions in place, this model explicitly captures and controls all the degradation root causes.
\end{itemize}

In our case study, we model degradation by using sensor data from a rotating machinery that is subjected to accelerated life testing (ALT) experiment \cite{gebraeel2005residual}. ALT experiments impose heightened loading/stress (e.g. acceleration factor) to observe changes in asset behaviour from brand new stage to failure. During these experiments, sensor readings, such as vibration, are continuously acquired, and are used to provide inferences on asset degradation. Degradation in rotating machinery (such as rolling element bearings) predominantly occur due to fatigue stresses. As the machinery progresses through its life, small cracks begin to materialize at the bearing raceway. These cracks progress and deteriorate over time until failure. This degradation process manifests itself through specific failure-induced vibration frequencies.  
The severity of these frequencies are used to form the degradation signal. When the degradation signal reaches the failure threshold, it induces a failure in the rotating machinery. In this case study, we used the degradation signal observations to derive parameters used for emulating degradation processes in our study.

Figure \ref{fig:synthetic_data_experiments_1} showcases the O\&M outcomes for different O\&M models. The plot showcases the impact of budget parameter on the value of the objective function for the robust optimization model, and the average simulation outcomes. It can be observed that increasing budget results in rising objective function value for the robust objective function. This is due to the additional costs of considering a wider set of degradation trajectories. As budget value increases, we observe that the simulation outcomes and the robust objective function value converges, showcasing that the proposed model effectively models and accounts for different simulation realizations with increasing budget values. Simulation outcomes obtain a minimum value at budget of $4$, after which the increasing level of conservatism starts to adversely impact the O\&M outcomes; resulting in unnecessary or early maintenances, and increased limitations on production values to account for degradation trajectories that are very unlikely to occur.

The table on right hand side of Figure \ref{fig:synthetic_data_experiments_1} showcases the performance of the deterministic benchmark models, compared to the deterministic version of our model. Compared to the benchmark models, the simulation outcomes for the proposed model provides significant benefits in O\&M costs owing to its accurate modeling of the underlying degradation process, and the capability to influence the degradation trajectory. Compared to Base, OID and MDI benchmark models, the deterministic counterpart of the proposed model provides $50.8\%$, $23.5\%$ and $41.1\%$ of O\&M cost improvements, respectively.

O\&M cost results from the proposed model arises as a function of improvements in operational and reliability metrics. In operations, Figure \ref{fig:synthetic_data_experiments_2} showcases that the unsatisfied demand due to asset outages are minimized for the proposed comprehensive O\&M model regardless of the budget. A parallel set of observations are also apparent in the number of failure instances \ref{fig:synthetic_data_experiments_2}. The proposed O\&M model reduces failure instances significantly, eventually converging to zero-failure outcomes at reasonable budget levels. The benchmark models incur failure instances that are significantly higher than the proposed model averages.

\subsection{Acceleration Method}

The proposed comprehensive O\&M model is a difficult problem to solve for large-scale instances. We develop an exact acceleration method to speed up the solution process to ensure that the proposed model can be solved for these realistic cases. Our approach revolves around generating upper bounds from ultra-conservative scenarios where all uncertain variables assume their maximum values as if there is no collective budget. We provide this integer solution as an efficient initial point. Details of the proposed method are outlined in the online appendix. 
Table \ref{acceleration_table} showcases experimental results that test the impact of the acceleration method in the solution performance for the comprehensive O\&M model. We cycle through a number of different parameter configurations that might impact solution performance to compare the method in different settings. The solution time results indicate that the proposed acceleration method improves number of solved instances and solution times.

\begin{table}[htbp]
  \centering
  {\tiny
  \caption{Comparative Results for the Performance of the Commercial Solver and the Acceleration Method}
    \begin{tabular}{c | c c c | c c c}
    
    \multicolumn{1}{c|}{} & \multicolumn{3}{c|}{\textbf{Commercial Solver}} & \multicolumn{3}{c}{\textbf{Acceleration Algorithm}} \\
    \hline
    \textbf{Budget} & \textbf{Avg. Time (s)} & \textbf{\# Solved} & \textbf{Avg. Gap (\%)} & \textbf{Avg. Time (s)} & \textbf{\# Solved} & \textbf{Avg. Gap (\%)} \\
    \hline
    \textbf{0.25} & 10535 & 6     & 1.10\% & 7861  & 8     & 0.66\% \\
    \hline
    \textbf{0.5} & 10595 & 7     & 1.14\% & 8225  & 8     & 0.79\% \\
    \hline
    \textbf{1} & 11643 & 4     & 1.39\% & 9179  & 8     & 1.10\% \\
    \hline
    \textbf{1.5} & 12711 & 3     & 1.57\% & 10741 & 5     & 1.25\% \\
    \hline
    \textbf{2} & 12801 & 4     & 1.63\% & 10230 & 6     & 1.15\% \\
    \hline
    \textbf{4} & 14378 & 1     & 2.21\% & 12904 & 4     & 1.54\% \\
    \hline
    \end{tabular}}
  \label{acceleration_table}
\end{table}%

\section{Conclusion}

In this article, we provide an integrated framework for O\&M of a fleet of assets by offering a robust optimization formulation that inherently captures degradation uncertainty through two modeling innovations: (i) reformulating and embedding continuous-time continuous-state asset degradation models and their connections to fleet-level performance as a set of linear constraints, and (ii) formulating a new generation of uncertainty sets for degradation that adapts to sensor-driven updates on asset degradation. The proposed model offers a significant deviation from the state-of-the-art O\&M models that either rely on static rules for maintenance (i.e., periodic) or do not provide computational scalability and multi-asset degradation interactions for large-scale systems.

The proposed O\&M policy achieves these improvements due to two factors. First, comprehensive modeling of the degradation process enables the O\&M model to accurately predict the failure risks while scheduling maintenance. Second, the proposed O\&M model also enables us to finetune the decision-dependent degradation factors to mitigate their impact on fleet-level O\&M. On the one hand, it adjusts operational decisions to control operational loading (i.e., OID effect) on asset degradation. On the other hand, it strategically schedules O\&M on connected assets to minimize and mitigate the degradation rate interactions across the fleet (i.e., MDI effect). The resulting policy yields a robust O\&M optimization model that synergizes operations and failure risks by accurately modeling degradation and operational interactions across a fleet of assets. It also showcases that sensor-driven comprehensive O\&M models can significantly impact operational outcomes.

The proposed model addresses a fundamental challenge in commercialization of condition-based maintenance systems: \textit{how to optimally translate asset-specific degradation insights to optimize fleet-level O\&M decisions for multiple interacting assets?} Most commercial applications and condition-based O\&M models focus on prediction of failure risks, and suggest simple asset-specific decisions: e.g., conduct maintenance when degradation exceeds a certain limit. These asset-specific decision policies perform poorly in complex multi-asset systems where O\&M actions require joint consideration of all the assets to achieve a fleet-level objective. Fleet-level O\&M models in literature typically use predetermined time windows to schedule maintenance, and do not take into account the degradation of the assets. By embedding degradation signals within a robust optimization formulation, we offer a new generation of large-scale O\&M models that can seamlessly translate degradation insights to fleet-level O\&M optimization, and close the gap across prediction and prescription.

The proposed optimization framework unlocks research directions to model a wide-range of fundamental problems in system reliability, and O\&M in complex systems; e.g., incorporation of multiple maintenance modes (e.g., minimal or partial maintenance), modeling of standby assets, and joint optimization of CBM and spare-part logistics for large-scale fleets.

\backmatter



\begin{appendices}

\section{Decision Dependent Degradation Factors}
\label{Decision Dependent Degradation Factors - Supplement}

In this section, we elucidate the decision-dependent degradation factors outlined in the main paper. We will first provide a detailed explanation and use cases for these degradation factors, followed by illustrations with a set of signals to showcase the impact of these factors on degradation and asset failure risks. As outlined in the main manuscript, we consider the following two types of decision-dependent degradation factors:

\underline{\textit{Operations Induced Degradation (OID)}} relates to the dynamic stress on assets as a function of operational decisions. The main premise of these models is that operational decisions, such as how much to produce or how fast to run, typically impose significant stress on the assets and increase the pace of degradation. These types of relationships are omnipresent in manufacturing and service settings \citep{uit2020condition}, where increasing production levels have a profound impact on the rate of degradation and time of failure in the assets. For instance, fatigue failure of gear systems exhibit this type of behavior across many industries, including but not limited to aerospace, automotive, and power systems. Studies have shown that the crack growth rate in gears evolves as a function of stress levels due to cyclic operational loading \citep{zhao2015integrated}.
Similar degradation behavior is also present in wind turbine blades. Studies using field measurements and simulation data to investigate the effects of turbulence intensity on the fatigue life of wind turbine blades, suggest that higher turbulence intensity results in faster blade degradation and reduced fatigue life \citep{bergami2014analysis}.  

\begin{figure}[h]
\centering
  \includegraphics[width=1\textwidth]{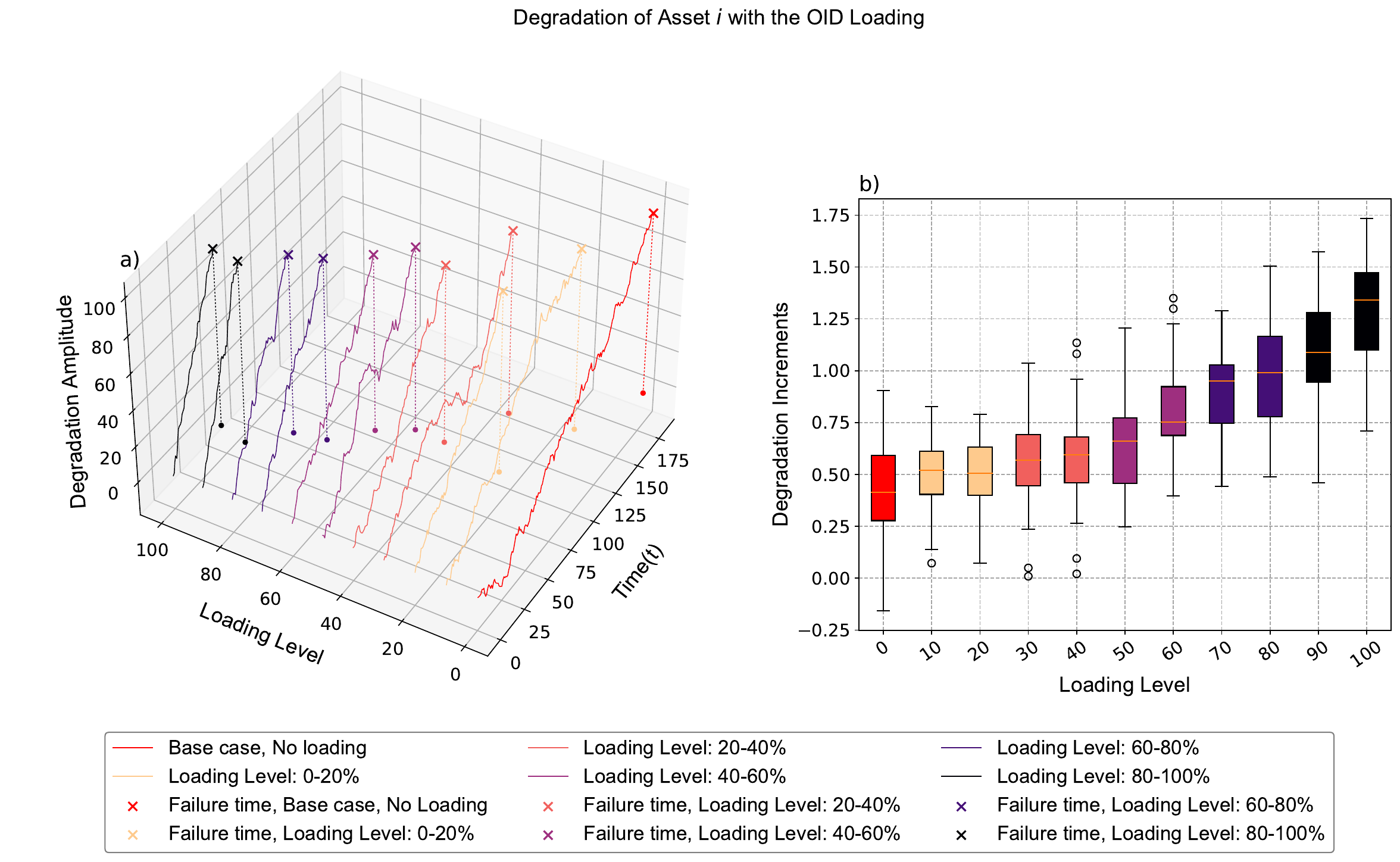}
  \caption{Impact of OID loading on asset degradation and failure risks. a) 3-D Plot for degradation signal realizations under different OID loading conditions, showcasing that increasing OID loading results in reductions in asset lifetime. b) Boxplot of degradation increments for different OID loading conditions, indicating that increasing loading leads to higher degradation increments and degradation rates.}
  \label{fig:OID_Plots}
\end{figure}

Figure \ref{fig:OID_Plots} presents two plots to showcase the behavior of degradation signals when subjected to different OID loading levels. Figure \ref{fig:OID_Plots}a presents the progression of the degradation signals when subjected to 11 different OID loading levels. Recall that a failure occurs when the degradation signal reaches a failure threshold, which corresponds to a degradation amplitude of 100. As expected, increasing loading levels mean that the assets degrade faster and reach the threshold earlier. Evidently, increasing loading reduces the lifetime of the assets. Changes in the rate (or pace) of degradation can be seen more clearly in Figure \ref{fig:OID_Plots}b. The plot showcases degradation increments, defined as the difference in degradation amplitudes across two degradation signal observations. This metric is an indicator of the rate of degradation. As expected, increasing OID loading causes a significant upward trend in degradation increments, and therefore, on the rate of degradation.

\textit{\underline{Multi-Asset Degradation Interactions (MDI)}} refers to degradation rate interactions across assets. We consider cases where an elevated degradation level of an asset sparks an increase in the rate of degradation of its connected assets. These type of interactions are omnipresent in any multi-asset (or multi-component) system with interactions. For instance, in wind turbine systems, degradation of hydrodynamic bearings leads to an increase in the looseness of primary transmission shafts. Consequently, this can raise the vibration levels in the gearbox and significantly impact the degradation of constituent gears  \citep{bian2014stochastic}. 
Another example comes from electro-mechanical cranking systems in cars. Cranking systems are the mechanisms that initiate the car engine upon starting the car (e.g. turning the key). Composed of a battery, starter motor, and an engine, cranking systems require all three components to function properly to achieve a successful cranking. When the starter motor degrades (e.g. brushes start to wear), batteries supply higher power and degrade faster to compensate for the degradation in the starter motor. A similar relationship also occurs in reverse. When the battery degrades, the starter motor operates longer to achieve a successful cranking. Degradation of the battery, and the starter motor, therefore become interdependent throughout their lifetimes.

\begin{figure}[h]
\centering
  \includegraphics[width=\textwidth]{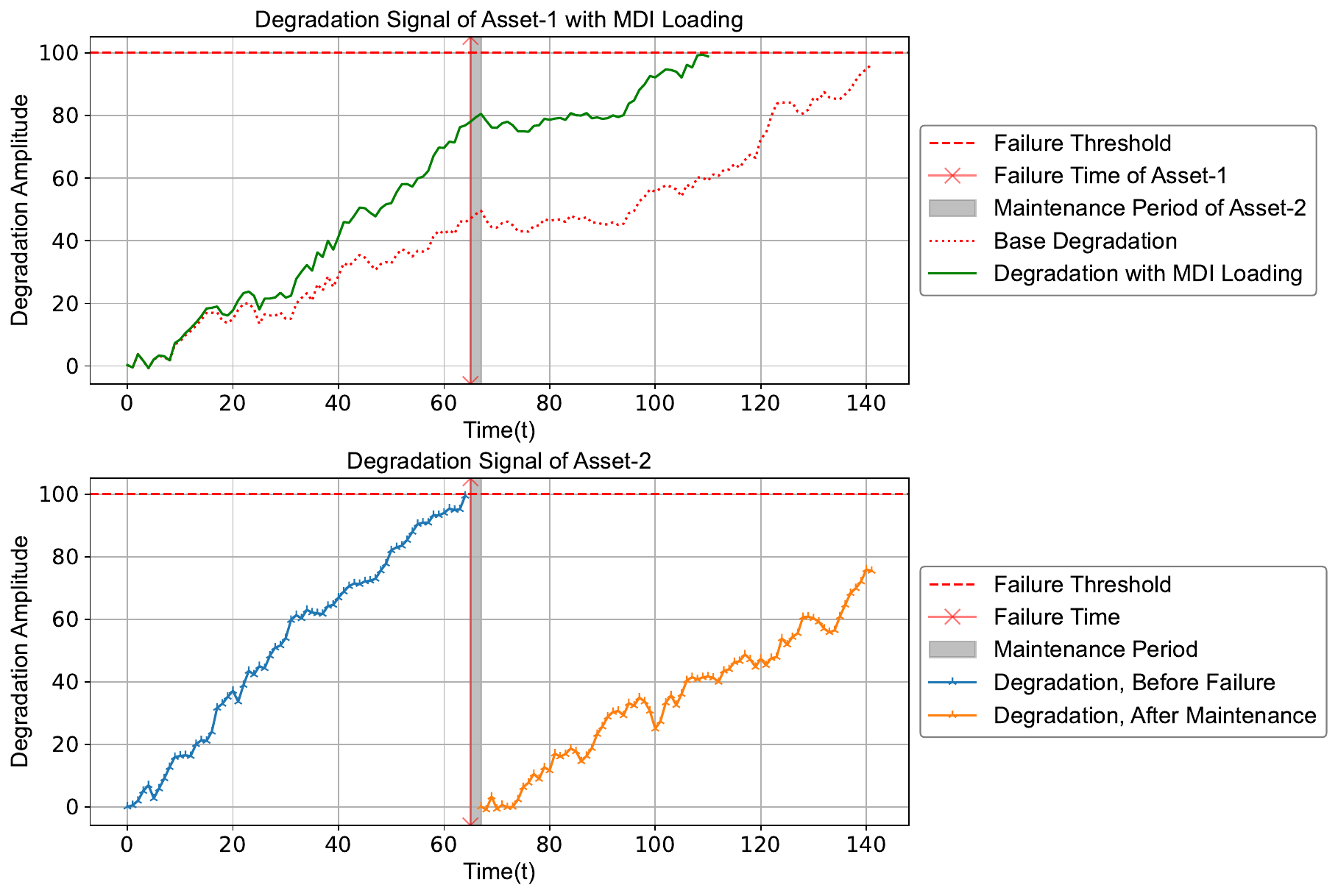}
  \caption{An Example on the impact of MDI loading on asset degradation and failure risks. Figure illustrates the progression of degradation signals for a case where asset-2, which undergoes maintenance at time $t=65$, impacts the rate of degradation in asset-1.}
  \label{fig:c2c}
\end{figure}

Figure \ref{fig:c2c} illustrates the MDI impact on degradation rates. The top and bottom plots of Figure \ref{fig:c2c}, shows how maintenance effect the signal amplitudes of asset-1 and asset-2 respectively. The gray line at $t = 65$ indicates the time of maintenance for asset-2. At the top plot, solid red line represents the observed degradation signal. The dashed line plots the base degradation, which is the degradation that would be observed if there was no MDI effect. As apparent in the figure, MDI impact caused a significant deviation in degradation across the base (no MDI) and observed degradation (with MDI), which results in significant changes to time of failure for asset-1. For the bottom plot, blue and orange signals indicate the degradation signal before and after maintenance, respectively. It can be observed that an increasing degradation signal amplitude in asset-2, resulted in higher rate of degradation in asset-1. In $t \in [40,65]$, asset-2 has higher degradation, hence resulting in a higher rate of degradation in asset 1. Conversely, in $t \in [67,90]$, asset-2 is in brand new condition, therefore the observed degradation rate in asset-1 is similar to its base degradation rate. 

\begin{figure}[h]
\centering
  \includegraphics[width=\textwidth]{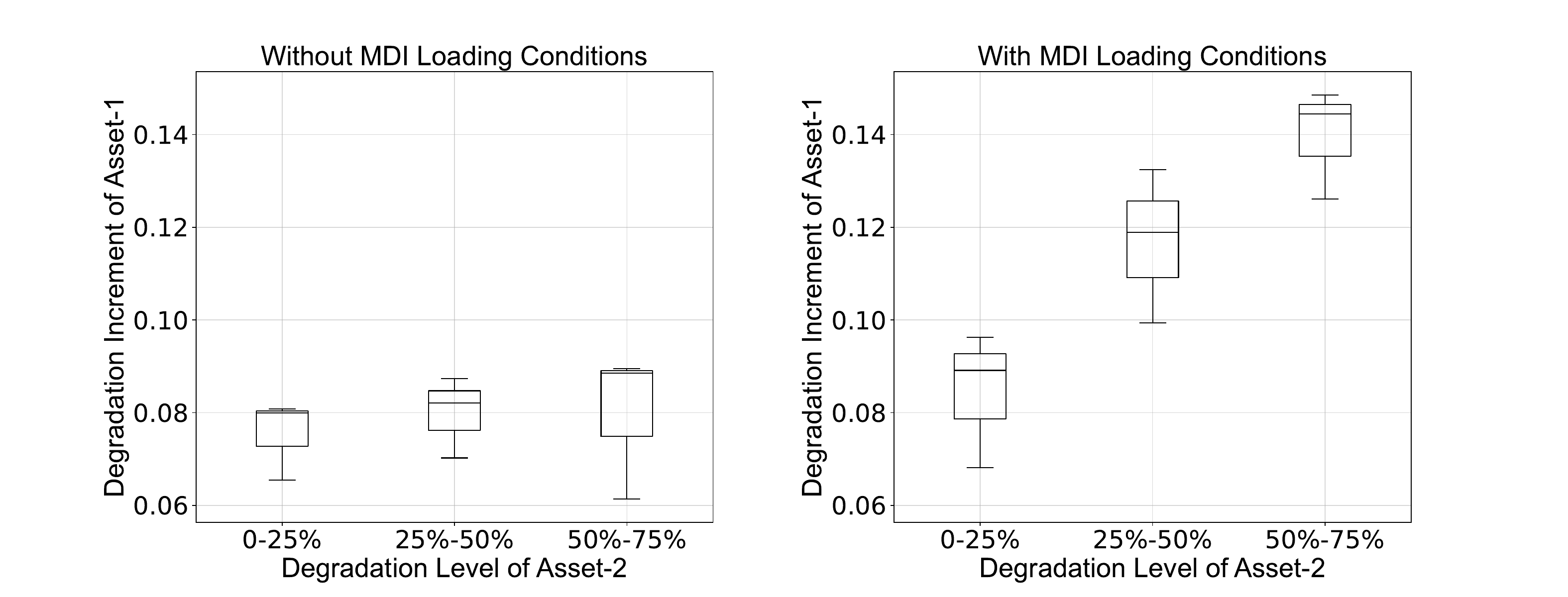}
  \caption{Impact of MDI loading on asset degradation and failure risks. a) Boxplot of degradation increments in asset-1 for different MDI loading conditions (i.e. percentile life of asset-2), where increasing loading leads to higher degradation increments and degradation rates. b) Corresponding degradation increments for the base degradation case showcasing no significant interaction between asset-1 and asset-2.}
  \label{fig:c2c_boxplot}
\end{figure}

To better illustrate this effect, Figure \ref{fig:c2c_boxplot} showcases the increments in degradation rate of asset-1, when asset-2 is in different percentiles of its useful life. The main observation is that when the life percentile of asset-2 increases (i.e. when asset-2 degrades), it should cause an aggravated rate of degradation in asset-1. This is in fact true. Figure \ref{fig:c2c_boxplot}a shows a boxplot of degradation increments in asset-1 as a function of lifetime percentile of asset-2. With increasing asset-2 life percentile, degradation increments exhibit an increasing trend. In contrast, Figure \ref{fig:c2c_boxplot}b showcases the corresponding results for the base case degradation model, where there is no significant correlation across changes in asset-1's degradation increment and asset-2's life percentile.

\setcounter{equation}{31}
\section{Sample Maintenance Schedule}
\begin{figure}[h!]
  \centering
  \includegraphics[width=0.8\textwidth]{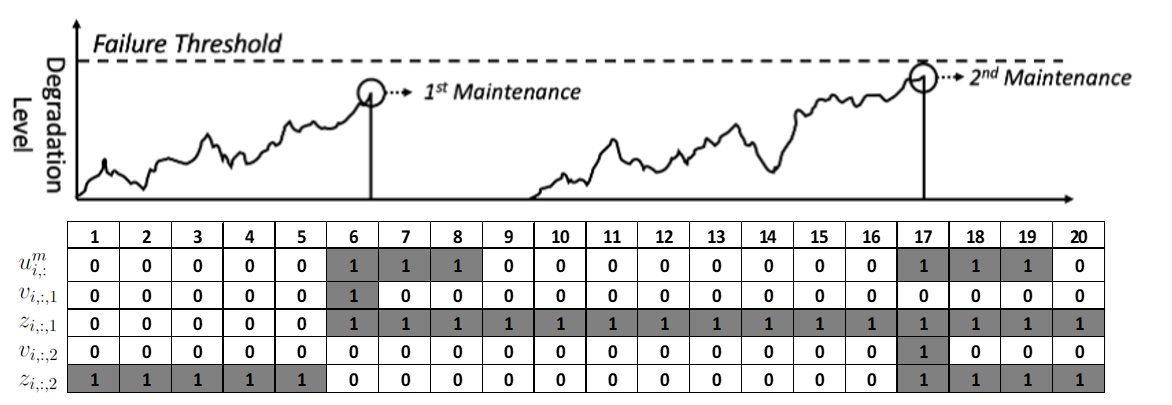}
  \caption{A sample maintenance schedule for Asset $i$ and its impact on the decision variables $v$, $z$ and $u^m$}
  \label{fig:example_maintenance}
\end{figure}

A significant improvement across the optimization models outlined in Sections 3.3. and 3.5. relates to the flexibility of the second model to allow scheduling multiple maintenance actions within a planning horizon. This change required us to revisit the decision variables to account for maintenance cycles. 

Figure \ref{fig:example_maintenance} demonstrates a sample configuration of maintenance decision variables $v$, $z$ and $u^m$. Upper graph shows how total degradation level evolves throughout the planning horizon. The table below shows corresponding maintenance decision variables  values. The variables $v_{i,:,1}$ and $v_{i,:,2}$ denote start times of the first and the second maintenance, respectively. The corresponding $z_{i,:,1}$ and $z_{i,:,2}$ are used to represent active maintenance cycles. Until the first maintenance occur, $z_{i,:,1}$ takes $0$, which means the first maintenance cycle is active. In period $6$, the first maintenance is initiated, $z_{i,:,1}$ start taking the value $1$. After the first maintenance, $z_{i,:,2}$ start taking the value $0$, which means the second maintenance cycle is active. This second maintenance cycle is active until the the second maintenance is scheduled in period $17$. The variable $u^m_{i,:}$ takes the value $1$ when the asset is unavailable due to a preventive maintenance actions (i.e. if it is undergoing maintenance), regardless of the corresponding maintenance cycle.

\section{Duality}
The proof of lemma 1 has three stages. In the first stage, we will develop an inner problem formulation for the following robust counterpart models of the degradation constraints:
\begin{align}
&l^{'}_{i,t,k}-l^{'}_{i,0,k}\geq \nonumber\\ &\hspace{0.7cm}\max_{\boldsymbol{d},\boldsymbol{\zeta},\boldsymbol{\gamma} \in \mathcal{U}_{i,t}} \left( \sum_{\tau = 1}^{t}d_{i,\tau}(1-z_{i,\tau,k}-u_{i,\tau}^{f})+  \sum_{\tau = 1}^{t}\zeta_{i,\tau}. p_{i,\tau,k}^{'}+ \sum_{\tau = 1}^{t}\sum_{j \in \mathcal{A}_i} \gamma_{j,i,\tau}.\omega_{j,i,\tau,k}^{'}\right), \nonumber&\\&\hspace{7.5cm} \forall i \in \mathcal{A},\hspace{2mm} \forall t \in \mathcal{T},\hspace{2mm}\forall k \in \mathcal{K}
\label{robust_degradation_supp}
\end{align}

Second stage of the proof uses the inner problem to generate a dual. Finally, in third stage, this dual formulation is incorporated into the rest of the optimization model. 

\underline{\textit{Stage 1:}} In right hand side of constraint \eqref{robust_degradation_supp}, we introduce the inner problem that maximizes degradation increments to search for the worst case realization within uncertainty set by controlling degradation parameters, $\boldsymbol{d}$, $\boldsymbol{\zeta}$ and $\boldsymbol{\gamma}$. Having transformed nested uncertainty sets into a set of constraints, we present the inner problem as follows:

\begin{subequations}
\begin{align}
\max \hspace{2mm}&\sum_{\tau = 1}^{t} \Bigg(d_{i,\tau}\cdot(1-z_{i,\tau,k}-u_{i,\tau}^{f})+ \zeta_{i,\tau}. p_{i,\tau}^{'}+\sum_{j \in \mathcal{A}_i} \gamma_{j,i,\tau}.\omega_{j,i,\tau,k}^{'}\Bigg)\\
s. t. & \nonumber\\
& \displaystyle \sum_{\tau = 1}^{t}\left( \left(\frac{1}{\hat{d}_{i,\tau}}\right) d_{i,\tau}+\left(\frac{1}{\hat{\zeta}_{i,\tau}}\right) \zeta_{i,\tau} +\sum_{j \in \mathcal{A}_i} (\frac{1}{\hat{\gamma}_{j,i,\tau}})\gamma_{j,i,\tau} \right) \leq  \Delta^{'}_{i,t}
\end{align}
\begin{align}
&d_{i,\tau} \leq \bar{d}_{i,\tau}+\hat{d}_{i,\tau}&\forall \tau \in [1,t]\\
&d_{i,\tau} \geq \bar{d}_{i,\tau} - \hat{d}_{i,\tau}&\forall \tau \in [1,t] \\
&\zeta_{i,\tau} \leq \bar{\zeta}_{i,\tau}+\hat{\zeta}_{i,\tau}&\forall \tau \in [1,t] \\
&\zeta_{i,\tau} \geq \bar{\zeta}_{i,\tau} - \hat{\zeta}_{i,\tau}&\forall \tau \in [1,t] \\
&\gamma_{j,i,\tau} \leq \bar{\gamma}_{j,i,\tau}+\hat{\gamma}_{j,i,\tau}& \forall j \in \mathcal{A}_i,\forall \tau \in [1,t]\\
&\gamma_{j,i,\tau} \geq \bar{\gamma}_{j,i,\tau} - \hat{\gamma}_{j,i,\tau}& \forall j \in \mathcal{A}_i,\forall \tau \in [1,t]
\end{align}\label{inner_problem}
\end{subequations}

where $\Delta_{i,t}^{'} = \Delta_{i,t} +\sum_{\tau = 1}^{t}\bigg( \frac{\bar{d}_{i,\tau}}{\hat{d}_{i,\tau}}+\frac{\bar{\zeta}_{i,\tau}}{\hat{\zeta}_{i,\tau}} +\sum_{j \in \mathcal{A}_i}\frac{\bar{\gamma}_{j,i,\tau}}{\hat{\gamma}_{j,i,\tau}}\bigg)$ for all $i \in \mathcal{A}$ and for all $t \in \mathcal{T}$. Note that, we control the level of conservatism by varying $\boldsymbol{\Delta}$ values: i.e. the bigger the $\boldsymbol{\Delta}$ values, the higher the level of conservatism. 

\underline{\textit{Stage 2:}} Integrating the inner problem \eqref{inner_problem} into main problem formulation enables us to hedge against possible disrupting realizations within the uncertainty set. To do so, we take the dual of the inner problem \eqref{inner_problem} to convert it into a minimization model, and 
build a compact problem formulation that jointly considers the exploration of worst-case degradation realizations and their associated impact on operations and maintenance decisions. We present the corresponding dual of the inner problem which replaces the of constraints \eqref{robust_degradation_supp}.  For all $i \in \mathcal{A}$, $t \in\mathcal{T}$, and $k \in \mathcal{K}$ the dual is formulated as follows:
\begin{subequations}\label{inner_obj}\begin{align}
\min \hspace{1mm} &\Delta^{'}_{i,t}\cdot(\pi^1_{i,t,k})+ \sum_{\tau=1}^{t}\Big[(\bar{d}_{i,\tau}+\hat{d}_{i,\tau})\cdot\pi^2_{i,t,\tau,k}-(\bar{d}_{i,\tau}-\hat{d}_{i,\tau}) \cdot\pi^3_{i,t,\tau,k}\nonumber\\&\hspace{1.5cm}+(\bar{\zeta}_{i,\tau}+\hat{\zeta}_{i,\tau})\cdot\pi^4_{i,t,\tau,k}-(\bar{\zeta}_{i,\tau}- \hat{\zeta}_{i,\tau})\cdot\pi^5_{i,t,\tau,k}\nonumber\\&\hspace{1.5cm}+\sum_{j \in \mathcal{A}_i}\big[(\bar{\gamma}_{j,i,\tau}+\hat{\gamma}_{j,i,\tau})\cdot\pi^6_{j,i,t,\tau,k}-(\bar{\gamma}_{j,i,\tau}-\hat{\gamma}_{j,i,\tau})\cdot\pi^7_{j,i,t,\tau,k}\big]\Big]\\
\text{s. t. } & \nonumber\\
&\frac{1}{\hat{d}_{i,t}}\pi^1_{i,t,k}+\pi^2_{i,t,\tau,k}-\pi^3_{i,t,\tau,k} \geq 1-z_{i,\tau,k}-u_{i,\tau}^{f},\hspace{0.5cm}\forall \tau \in [1,t]\label{71}\\
&\frac{1}{\hat{\zeta}_{i,t}}\pi^1_{i,t,k}+\pi^4_{i,t,\tau,k}-\pi^5_{i,t,\tau,k} \geq p^{'}_{i,\tau,k},\hspace{2.2cm}\forall \tau \in [1,t]\label{72}\\
&\frac{1}{\hat{\gamma}_{j,i,\tau}}\pi^1_{i,t,k}+\pi^6_{j,i,t,\tau,k}-\pi^7_{j,i,t,\tau,k} \geq\omega^{'}_{j,i,\tau,k},\hspace{1.3cm}\forall j \in \mathcal{A}_i,\forall \tau \in [1,t]\label{73}\\
&\pi^1_{i,t,\tau,k},\pi^2_{i,t,\tau,k},\pi^3_{i,t,\tau,k},\pi^4_{i,t,\tau,k},\pi^5_{i,t,\tau,k},\pi^6_{j,i,t,\tau,k},\pi^7_{j,i,t,\tau,k} \in {\rm I\!R^{+}}\label{74}
\end{align}
\end{subequations}
The resulting formulation investigates worst-case realizations within the uncertainty set so that we can immunize our decisions against uncertainty.

\underline{\textit{Stage 3:}} Finally, we incorporate the dual formulation within the robust optimization model. The resulting optimization model can be cast as follows:
\newpage

\begin{align}
    \min \quad & \eqref{objective function} \nonumber\\
    \text{s.t.}\nonumber\\ \quad & l^{'}_{i,t,k}-l^{'}_{i,0,k} \geq \left\{ \begin{aligned}
        &\Delta^{'}_{i,t}\cdot(\pi^1_{i,t,k}) \\        &+\sum_{\tau=1}^{t}\Big[(\bar{d}_{i,\tau}+\hat{d}_{i,\tau})\cdot\pi^2_{i,t,\tau,k}-(\bar{d}_{i,\tau}-\hat{d}_{i,\tau}) \cdot\pi^3_{i,t,\tau,k}\\        &+(\bar{\zeta}_{i,\tau}+\hat{\zeta}_{i,\tau})\cdot\pi^4_{i,t,\tau,k} -(\bar{\zeta}_{i,\tau}-\hat{\zeta}_{i,\tau})\cdot\pi^5_{i,t,\tau,k}\\        &+\sum_{j \in \mathcal{A}_i}\big[(\bar{\gamma}_{j,i,\tau}+\hat{\gamma}_{j,i,\tau})\cdot\pi^6_{j,i,t,\tau,k}\\        &\qquad-(\bar{\gamma}_{j,i,\tau}-\hat{\gamma}_{j,i,\tau})\cdot\pi^7_{j,i,t,\tau,k}\big]\Big]
\end{aligned}\right\},\label{dual_degradation_update123}\nonumber\\&\hspace{7.2cm}  \forall i \in \mathcal{A}, t \in \mathcal{T}, k \in \mathcal{K}  \\
    & (\ref{component_2}-\ref{18})\nonumber\\&(\ref{threshold2}-\ref{10})\nonumber\\&(\ref{71})-(\ref{74})\nonumber
\end{align}\label{Robust Counterpart}
This formulation can be written in a condensed form as follows:
\begin{align}
    Min \{ \;\eqref{objective function} \; \nonumber \\  \; s.t. \nonumber\\&\;  (\ref{component_2}-\ref{18}),\nonumber\\&(\ref{threshold2}-\ref{10}) \nonumber, \\&\mathcal{L} \geq \boldsymbol{\pi} \mathbf{u} ,\boldsymbol{\pi}\mathbf{d} \geq \mathbf{m} \}
\end{align}
where $\mathcal{L}$ denotes left hand side of constraints 
(35), {$\boldsymbol{\pi}$} includes corresponding dual variables ($\pi^1_{i,t,k}, \pi^2_{i,t,\tau,k}, \pi^3_{i,t,\tau,k}, \pi^4_{i,t,\tau,k}, \pi^5_{i,t,\tau,k},\pi^6_{j,i,t,\tau,k},\pi^7_{j,i,t,\tau,k}$) for each constraint in uncertainty set (28), $\mathbf{u}$ represents the uncertainty sets parameters $\Delta^{'}_{i,t},(\bar{d}_{i,\tau}+\hat{d}_{i,\tau}),(\bar{d}_{i,\tau}-\hat{d}_{i,\tau)} ,(\bar{\zeta}_{i,\tau}+\hat{\zeta}_{i,\tau}),(\bar{\zeta}_{i,\tau}-\hat{\zeta}_{i,\tau}),(\bar{\gamma}_{j,i,\tau}+\hat{\gamma}_{j,i,\tau}),(\bar{\gamma}_{j,i,\tau}-\hat{\gamma}_{j,i,\tau})$ which appears in the inner minimization problem objective function \eqref{inner_obj} and the right hand side of constraints (35).
$\mathbf{d}$ represents coefficient matrix of constraints (\ref{71}-\ref{73}). $\mathbf{m}$ represents the maintenance and degradation variables which constitute right hand side coefficients of constraints (\ref{71}-\ref{73}).

\section{Acceleration Method}
The comprehensive O\&M model poses a challenge when attempting to solve it for large-scale problem instances. We employ an acceleration method to speed up the solution process, while also maintaining the optimality guarantees.
The acceleration method outlined in Algorithm \ref{alg:accelerate} involves two initial steps which work in parallel. In step one, we provide an initial feasaible integer solution to the solver to start it from an efficient initial point. To achieve this, we first create a scenario where its deterministic solution would yield the upper bound to the optimal solution. We choose uncertain degradation parameters $\boldsymbol{d}$, $\boldsymbol{\zeta}$, and $\boldsymbol{\gamma}$ at their maximum value (by assuming that the budget  is infinity) to ensure that the resulting scenario solution would be an upper bound for the optimal solution. In other words, this solution will be feasible under all possible scenarios. 
Having solved the deterministic problem of this ultra-conservative scenario, we provide an integer solution to the solver. This helps solver to reach an optimal solution faster. In step two, we choose a subset of realizations for the degradation parameters $\boldsymbol{d}$, $\boldsymbol{\zeta}$, and $\boldsymbol{\gamma}$ considering budget $\boldsymbol{\Delta}$ so that resulting scenario would be one of the candidate worst-case scenarios that model needs to immunize against. We add a number of these extreme scenarios in the form of constraints (5). Including such scenarios in the robust model improve solver performances. Once we provide an integer solution and add extreme scenarios to the robust model, we solve the resulting optimization problem. Implementation is built using  Gurobi 9.0.3 and Python with Intel(R) Core(TM) i7-9850H CPU @ 2.60GHz with 16 GB RAM. A pseudocode associated with the proposed solution algorithm is presented in Algorithm \ref{alg:accelerate}.

\begin{algorithm}
	\algtext*{Indent}
	\caption{ Acceleration Method}\label{alg:accelerate}
	\begin{algorithmic}[]
		\State \textbf{Step 1: Solve Ultra Conservative Deterministic Problem to Get an Initial Solution}
		\State  \indent Generate a scenario where all uncertain variables assume their maximum values
		\State  \indent Solve deterministic optimization problem using the extreme scenario
		\State  \indent Obtain integer variable solutions		
		\State \textbf{Step 2: Add a Subset of Worst-Case Scenario Cuts to Accelerate the Convergence}
		\State  \indent Generate extreme scenarios
		\State  \indent Add cuts of these scenarios in form of the degradation update constraints 
		\State \textbf{Step 3: Solve the Resulting Robust Optimization Problem}
		\State \indent Provide the integer solution from Step 1 as warm-start values
		\State \indent Solve robust optimization problem

	\end{algorithmic}
\end{algorithm}

\end{appendices}


\bibliography{sn-article}

\end{document}